\documentclass[11pt,reqno]{article}
\usepackage{amsmath, amssymb}
\numberwithin{equation}{section}

\newtheorem{theorem}{Theorem}[section]
\newtheorem{definition}[theorem]{Definition}
\newtheorem{proposition}[theorem]{Proposition}
\newtheorem{corollary}[theorem]{Corollary}
\newtheorem{lemma}[theorem]{Lemma}
\newtheorem{conjecture}[theorem]{Conjecture}
\newtheorem{fact}[theorem]{Fact}

\def \proof {\noindent {\bf Proof.}\ \ }
\def \remark {\noindent {\bf Remark.}\ \ }
\def \remarks {\noindent {\bf Remarks.}\ \ }

\def \endproof {{\mbox{}\nolinebreak\hfill\rule{2mm}{2mm}\par\medbreak}}

\def \N {\mathbb{N}}
\def \R {\mathbb{R}}

\def \Z {\mathbb{Z}}
\def \E {\mathbb{E}}
\def \P {\mathbb{P}}

\def \one {{\bf 1}}

\def \CC {\mathcal{C}}

\def \OO {\mathcal{O}}
\def \PP {\mathcal{P}}

\def \QQ {\mathcal{Q}}
\def \F  {{F}}
\def \a {\alpha}
\def \b {\beta}
\def \g {\gamma}
\def \e {\varepsilon}

\def \d {\delta}

\def \l {\lambda}
\def \L {\Lambda}
\def \s {\sigma}

\def \w {\omega}
\def \W {\Omega}
\def \< {\langle}
\def \> {\rangle}

\def \sign {{\rm sign}}

\def \rank {{\rm rank }}

\def \cconv {{\rm cconv}}

\def \id {{\it id}}

\def \vol {{\rm vol}}

\def \vc {{\it v}}
\def \tower {{\rm Tower}}
\def \Ball {{\rm Ball}}
\def \rad {{\it rad}}

\begin{document}
\title {Combinatorics of random processes
        and sections \\ of convex bodies}
\author {M. Rudelson\footnote{
   Department of Mathematics,
   University of Missouri,
   Columbia, MO 65211, USA;
   e-mail: \mbox{rudelson@math.missouri.edu} }
   \and
   R. Vershynin\footnote{
   Deptartment of Mathematics,
   University of California,
   Davis, CA 95616, USA;
   e-mail: \mbox{vershynin@math.ucdavis.edu}} }

\maketitle

\begin{abstract}
We find a sharp combinatorial bound for the metric entropy of sets in
$\R^n$ and general classes of functions. This solves two basic
combinatorial conjectures on the empirical processes. 1. A class of
functions satisfies the uniform Central Limit Theorem if the
square root of its combinatorial dimension is integrable. 2. The
uniform entropy is equivalent to the combinatorial dimension under
minimal regularity. Our method also constructs a nicely bounded
coordinate section of a symmetric convex body in $\R^n$. In the
operator theory, this essentially proves for all normed spaces the
restricted invertibility principle of Bourgain and Tzafriri.
\end{abstract}




\qquad

\section{Introduction}                    \label{introduction}

This paper develops a sharp combinatorial method
for estimating metric entropy of sets in $\R^n$
and, equivalently, of function classes on a probability space.
A need in such estimates occurs naturally in a
number of problems of analysis (functional, harmonic and
approximation theory), probability, combinatorics, convex and
discrete geometry, statistical learning theory, etc.
Our entropy method, which evolved from the work of 
S.Mendelson and the second author \cite{MV},
is motivated by several problems in the empirical 
processes, asymptotic convex geometry and operator theory.

Throughout the paper, $F$ is a class of real valued functions
on some domain $\Omega$. It is a central problem of the theory
of empirical processes to determine
whether the classical limit theorems hold uniformly over $F$.
Let $\mu$ be a probability distribution on $\Omega$ and
$X_1, X_2, \ldots \in \Omega$ be independent samples
distributed according to a common law $\mu$.
The problem is to determine whether the sequence
of real valued random variables $(f(X_i))$ obeys the
central limit theorem uniformly over all $f \in F$ and
over all underlying probability distributions $\mu$, i.e. whether
the random variable
$\frac{1}{\sqrt{n}} \sum_{i=1}^n (f(X_i) - f(X_1))$ converges
to a Gaussian random variable uniformly. With the right definition of
the convergence, if that happens, $F$ is a {\em uniform Donsker class}.
The precise definition can be found in \cite{LT} and \cite{Du 99}.

The pioneering work of Vapnik and Chervonenkis
\cite{VC 68, VC 71, VC 81} demonstrated that the validity
of the uniform limit theorems on $F$ is connected with
the combinatorial structure of $F$, which is quantified by what we
call the {\em combinatorial dimension} of $F$.
For classes of $\{0,1\}$-valued functions,
it is the classical Vapnik-Chervonenkis dimension.
For a general class $F$ and $t \ge 0$,
a subset $\s$ of $\Omega$ is called $t$-shattered
by a class $F$ if there exists a level function $h$ on $\s$
such that, given any partition $\s = \s_- \cup \s_+$,
one can find a function $f \in F$ with $f(x)  \le  h(x)$ if
$x \in \s_-$ and $f(x)  \ge  h(x) + t$ if $x \in \s_+$.
The combinatorial dimension of $F$, denoted by $\vc(F,t)$,
is the maximal cardinality of a set $t$-shattered by $F$.
Simply speaking, $\vc(F,t)$ is the maximal size of a set
on which $F$ oscillates in all possible $\pm t/2$ ways
around some level $h$.

Connections between the combinatorial dimension (and its variants)
with the limit theorems of probability theory have been the major 
theme of many papers. For a comprehensive account of what was 
known about these profound connections by 1999, 
we refer the reader who to the book of Dudley \cite{Du 99}.

Dudley proved that a class $F$ of $\{0,1\}$-valued functions
is a uniform Donsker class if and only if its combinatorial
(Vapnik-Chernovenkis) dimension $\vc(F,1)$ is finite. 
This is one of the main results on the empirical processes 
for $\{0,1\}$ classes. The problem for general classes turned
out to be much harder \cite{T 02}, \cite{MV}.
In the present paper we prove an optimal integral description
of uniform Donsker classes in terms of the combinatorial dimension.

\begin{theorem}                             \label{i:donsker}
  Let $F$ be a uniformly bounded class of functions.
  Then
  $$
  \int_0^\infty \sqrt{\vc(F,t)} \; dt < \infty
  \ \Rightarrow \ \text{$F$ is uniform Donsker}
  \ \Rightarrow \ \vc(F,t) = O(t^{-2}).
  $$
\end{theorem}

This trivially contains Dudley's theorem on the $\{0,1\}$ classes.
M.Talagrand proved Theorem \ref{i:donsker} with an extra factor of
$\log^M(1/t)$ in the integrand and asked about the optimal value
of the absolute constant exponent $M$ \cite{T 92}, \cite{T 02}.
Talagrand's proof was based on a very involved iterational argument. 
In \cite{MV}, S.Mendelson and the second author introduced
a new combinatorial idea. Their approach led to a much clearer
proof, which allowed to reduce the exponent to $M = 1/2$. 
Theorem \ref{i:donsker} removes the logarithmic factor completely, 
thus the optimal exponent is $M = 0$. Our argument significantly
relies on the ideas originated in \cite{MV} and also uses 
a new iterational method.   
The second implication of Theorem \ref{i:donsker}, which makes sense 
for $t \to 0$, is well-known (\cite{Du 99} 10.1).

Theorem \ref{i:donsker} reduces to estimating metric entropy of $F$
by the combinatorial dimension of $F$.
For $t > 0$, the {\em Koltchinskii-Pollard entropy} of $F$ is
\[
  D(F, t)= \log \sup \Big( n \mid \exists f_1, \ldots, f_n \in F \ \
  \forall i < j \ \int (f_i-f_j)^2 d \mu \ge t^2 \Big)
\]
where the supremum is by $n$ and over all probability measures
$\mu$ supported by the finite subsets of $\W$. It is easily seen
that $D(F,t)$ dominates the combinatorial dimension: $D(F,t)
\gtrsim \vc(F,2t)$. Theorem \ref{i:donsker} should then be
compared to the fundamental description valid for all uniformly bounded
classes:
\begin{equation}                               \label{i:dudley donsker}
\int_0^\infty \sqrt{D(F,t)} \; dt < \infty
\ \Rightarrow \ \text{$F$ is uniform Donsker}
\ \Rightarrow \ \ D(F,t) = O(t^{-2}).
\end{equation}
The left part of \eqref{i:dudley donsker} is a strengthening
of Pollard's central limit theorem and is due to Gine and Zinn
(see \cite{GZ}, \cite{Du 99} 10.3, 10.1). The right part is an observation
due to Dudley (\cite{Du 99} 10.1).

An advantage of the combinatorial description in Theorem \ref{i:donsker}
over the entropic description in \eqref{i:dudley donsker}
is that the combinatorial dimension is much easier to bound than
Koltchinskii-Pollard entropy (see \cite{AB}).
Large sets on which $F$ oscillates in all $\pm t/2$ ways are so
sound structures that their existence can be hopefully easily
detected or eliminated, which leads to an estimate on
the combinatorial dimension.
In contrast to this, bounding Koltchinskii-Pollard entropy involves
eliminating all large separated configurations $f_1,\ldots,f_n$
with respect to all probability measures $\mu$; this can be a hard
problem even on the plane (for a two-point domain $\Omega$).

The nontrivial part of Theorem \ref{i:donsker}
follows from \eqref{i:dudley donsker}
and the central result of this paper:

\begin{theorem}                                 \label{i:integral=integral}
  For every class $F$,
  \[
  \int_0^\infty \sqrt{D(F,t)} \; dt
  \asymp \int_0^\infty \sqrt{\vc(F,t)} \; dt.
  \]
\end{theorem}
The equivalence $\asymp$ is up to an absolute constant factor $C$,
thus $a \asymp b$ iff $a/C \le b \le C a$.

Looking at Theorem \ref{i:integral=integral} one
naturally asks whether the Koltchinskii-Pollard entropy is poinwise
equivalent to the combinatorial dimension.
M.~Talagrand indeed proved this for uniformly bounded classes
under minimal regularity and up to a logarithmic factor.
For the moment, we consider a simpler version of this
regularity assumption: there exists an $a > 1$ such that
\begin{equation}                        \label{i:min reg}
\vc(F,at) \le \frac{1}{2} \, \vc(F,t)
\ \ \ \text{for all $t > 0$.}
\end{equation}
In 1992, M.~Talagrand proved essentially under \eqref{i:min reg}
that for $0 < t < 1/2$
\begin{equation}                        \label{i:talagrand}
c \; \vc(F,2t) \le D(F,t) \le C \; \vc(F,ct) \log^{M}(1/t)
\end{equation}
\cite{T 92}, see \cite{T 87}, \cite{T 02}. Here $c>0$ is an
absolute constant and $M$ depends only on $a$. The question on the
value of the exponent $M$ has been open. S.~Mendelson and the
second author proved \eqref{i:talagrand} {\em without} the minimal
regularity assumption \eqref{i:min reg} and with $M = 1$, which is
an optimal exponent in that case. The present paper proves that
{\em with} the minimal regularity assumption, the exponent reduces
to $M = 0$, thus completely removing both the boundedness
assumption and the logarithmic factor from Talagrand's inequality
\eqref{i:talagrand}. As far as we know, this unexpected fact was
not even conjectured.

\begin{theorem}                         \label{i:entropy=dimension}
  Let $F$ be a class which satisfies the minimal regularity
  assumption \eqref{i:min reg}. Then for all $t > 0$
  $$
  c \; \vc(F,2t) \le D(F,t) \le C \; \vc(F,ct),
  $$
  where $c>0$ is an absolute constant and $C$
  depends only on $a$ in \eqref{i:min reg}.
\end{theorem}
Therefore, in presence of minimal regularity,
the Koltchinski-Pollard entropy and the combinatorial dimension are
equivalent. Rephrasing M.Talagrand's comments from \cite{T 02}
on his inequality \eqref{i:talagrand}, Theorem \ref{i:entropy=dimension}
is of the type ``concentration of pathology''.
Suppose we know that $D(F,t)$ is large. This simply means that
$F$ contains many well separated functions,
but we know very little about what kind of pattern
they form. The content of Theorem \ref{i:entropy=dimension}
is that it is possible to construct
a large set $\sigma$ on which not only many functions in $F$
are well separated from each other,
but on which they oscillate in {\em all} possible
$\pm ct$ ways. We now have a very precise structure that
witnesses that $F$ is large. This result is exactly in the line
of Talagrand's celebrated characterization of Glivenko-Cantelli
classes \cite{T 87}, \cite{T 96}.

Theorem \ref{i:entropy=dimension} remains true
if one replaces the $L_2$ norm
in the definition of the Koltchinski-Pollard entropy
by the $L_p$ norm for $1 \le p < \infty$.
The extremal case $p = \infty$ is important and more difficult.
The $L_\infty$ entropy is naturally
\[
  D_\infty(F, t)= \log \sup \Big( n \mid \exists f_1, \ldots, f_n \in F \ \
  \forall i < j  \ \sup_
\w |(f_i-f_j)(\w)| \ge t \Big).
\]
Assume that $F$ is uniformly bounded (in absolute value) by $1$.
Even then $D_\infty(F,t)$ can not be bounded by a function of $t$ and
$\vc(F, ct)$: to see this, it is enough to take for $F$ the
collection of the indicator functions of the intervals $[2^{-k-1},
2^{-k}]$, $k \in \N$, in $\Omega = [0,1]$. However, if $\Omega$ is
finite, it is an open question how the $L_\infty$ entropy depends
on the size of $\Omega$. N.Alon et al. \cite{ABCH} proved that if
$|\Omega| = n$ then $D_\infty(F, t) = O(\log^2 n)$ for fixed $t$
and $\vc(F, ct)$. They asked whether the exponent $2$ can be
reduced. We answer this by reducing $2$ to any number larger
than  the minimal possible value $1$. For every $\e \in (0,1)$,
\begin{equation}                                \label{i:by linfty}
D_\infty(F, t) \le  C v \, \log(n/vt) \cdot \log^\e(n/v),
\ \ \ \text{where $v = \vc(F, c \e t)$}
\end{equation}
and where $C, c > 0$ are absolute constants.
One can look at this estimate as a continuous asymptotic version
of Sauer-Shelah Lemma. The dependence on $t$ is optimal,
but conjecturally the factor $\log^\e(n/v)$ can be removed.

The combinatorial method of this paper applies to the study of coordinate
sections of a symmetric convex body $K$ in $\R^n$. The average
size of $K$ is commonly measured by the so-called M-estimate,
which is $M_K = \int_{S^{n-1}} \|x\|_K \; d \s(x)$, where $\s$ is
the normalized Lebesgue measure on the unit Euclidean sphere
$S^{n-1}$ and $\|\cdot\|_K$ is Minkowski functional of $K$.
Passing from the average on the sphere to the Gaussian average on
$\R^n$, Dudley's entropy integral connects the M-estimate to the integral
of the metric entropy of $K$; then Theorem
\ref{i:integral=integral} replaces the entropy by the
combinatorial dimension of $K$. The latter has a remarkable
geometric representation, which leads to the following result. For
$1 \le p \le \infty$ denote by $B_p^n$ the unit ball of the space
$\ell_p^n$:
$$
B_p^n = \{ x \in \R^n \;:\; |x_1|^p + \cdots + |x_n|^p \le 1 \}.
$$
If $M_K$ is large (and thus $K$ is small ``in average'') then
there exists a coordinate section of $K$ contained in the
normalized octahedron $D = \sqrt{n} B_1^n$. Note that the $M_D$ is
bounded by an absolute constant. In the rest of the paper, $C, C',
C_1, c, c', c_1, \ldots$ will denote positive absolute constants
whose values may change from line to line.

\begin{theorem}                                  \label{dvoretzky}
  Let $K$ be a symmetric convex body containing the unit
  Euclidean ball $B_2^n$, and let $M = c M_K \log^{-3/2}(2/M_K)$.
  Then there exists a subset $\s$ of $\{1, \ldots, n\}$
  of size $|\s|  \ge  M^2 n$, and such that
  \begin{equation}                               \label{inclusion}
    M \, (K \cap \R^\s) \subseteq \sqrt{|\s|} B_1^\s.
  \end{equation}
\end{theorem}
Recall that the classical Dvoretzky theorem in the form of Milman
guarantees, for $M = M_K$, the existence of a subspace $E$
of dimension $\dim E  \ge  c M^2 n$ and such that
\begin{equation}                                \label{i:dvoretzky}
c_1 B_2^n \cap E
\subseteq M (K \cap E)
\subseteq c_2 B_2^n \cap E.
\end{equation}
To compare the second inclusion of \eqref{i:dvoretzky}
to \eqref{inclusion}, recall that by Kashin's theorem (\cite{K 77},
\cite{K 85}, see \cite{Pi} 6) there exists a subspace $E$ in $\R^\s$
of dimension at least $|\s|/2$ such that
the section $\sqrt{|\s|} B_1^\s \cap E$ is equivalent to $B_2^n \cap E$.

A reformulation of Theorem \ref{dvoretzky} in the operator
language generalizes the restricted invertibility principle of
Bourgain and Tzafriri \cite{BT 87} to all normed spaces. Consider
a linear operator $T : l_2^n \to X$ acting from the Hilbert space
into arbitrary Banach space $X$. The ``average'' largeness of such
an operator is measured by its $\ell$-norm, defined as $\ell(T)^2
= \E \|T g\|^2$, where $g = (g_1, \ldots, g_n)$ and $g_i$ are
normalized independent Gaussian random variables. We prove that if
$\ell(T)$ is large then $T$ is well invertible on some large
coordinate subspace. For simplicity, we state this here for spaces
of type $2$ (see \cite{LT} 9.2), which includes for example
all the $L_p$ spaces and their subspaces for $2 \le p < \infty$. For
general spaces, see Section \ref{s:sections}.

\begin{theorem}[General Restricted Invertibility]  \label{i:invertibility}
  Let $T : l_2^n \to X$ be a linear operator with $\ell(T)^2 \ge n$,
  where $X$ is a normed space of type $2$.
  Let $\a = c \log^{-3/2}(2\|T\|)$.
  Then there exists a subset $\sigma$ of $\{1, \ldots, n\}$
  of size $|\s| \ge \a^2 n / \|T\|^2$ and such that
  $$
  \|Tx\|  \ge  \a \b_X \|x\| \ \ \
  \text{for all $x \in \R^\s$}
  $$
  where $c>0$ is an absolute constant and $\b_X > 0$ depends
  on the type $2$ constant of $X$ only.
\end{theorem}

Bourgain and Tzafriri essentially proved this restricted
invertibility principle for $X = l_2^n$
(and without the logarithmic factor),
in which case $\ell(T)$ equals the Hilbert-Schmidt norm of $T$.

The heart of our method is a result of combinatorial geometric
flavor. We compare the covering number of a convex body $K$ by a
given convex body $D$ to the number of the integer cells contained
in $K$ and its projections. This will be explained in detail in
Section~\ref{s:overview}. All main results of this paper are then
deduced from this principle. The basic covering result of this
type and its proof occupies Section~\ref{s:by tower}. First
applications to covering $K$ by ellipsoids and cubes appear in
Section \ref{s:by ellipsoids and cubes}. Estimate \eqref{i:by
linfty} is also proved there. Section~\ref{s:by lorentz} deals
with covering by balls of a general Lorentz space; the
combinatorial dimension controls such coverings. From this we
deduce in Section~\ref{s:gaussian} our main results, Theorems
\ref{i:integral=integral} and \ref{i:entropy=dimension}.
Theorem \ref{i:integral=integral} shows in particular that
{\em in the classical Dudley's entropy integral, the entropy can
be replaced by the combinatorial dimension}. This yields a
new powerful bound on Gaussian processes (see Theorem \ref{supremum}
below), which is a quantitative version of Theorem \ref{i:donsker}.
This method is used in Section~\ref{s:sections} to prove
Theorem~\ref{dvoretzky} on the coordinate sections of convex bodies.
Theorem~\ref{dvoretzky} is equivalently expressed in the operator
language as a general principle of restricted invertibility,
which implies Theorem~\ref{i:invertibility}.

\qquad

ACKNOWLEDGEMENTS. This project started when both authors
visited the Pacific Institute of Mathematical Sciences. We would
like to thank PIMS for its hospitality.
A significant part of the work was done when the second author was
PIMS Postdoctoral Fellow at the University of Alberta.
He thanks this institution and especially
Nicole Tomczak-Jaegermann for support and encouragement.

\section{The Method}                \label{s:overview}

Let $K$ and $D$ be convex bodies in $\R^n$. We are interested in
the covering number $N(K,D)$, the minimal number of translates of
$D$ needed to cover $K$. More precisely, $N(K,D)$ is the minimal
number $N$ for which there exist points $x_1, x_2, \ldots x_N$
satisfying
$$
K \subseteq \bigcup \limits_{j=1}^N (x_j+D).
$$
Computing the covering number is a very difficult problem
even in the plane \cite{CFG}. Our main idea is to relate the
covering number to the {\em cell content} of $K$, which we define
as the number of the the integer cells contained in all coordinate
projections of $K$:
\begin{equation}                            \label{lattice content}
\Sigma(K) = \sum_P \text{number of integer cells contained in
$PK$}.
\end{equation}
The sum is over all $2^n$ coordinate projections in $\R^n$, i.e.
over the orthogonal projections $P$ onto $\R^\s$ with $\s \subseteq
\{1,\ldots, n\}$. The integer cells are the unit cubes with integer
vertices, i.e. the sets of the form
$a + [0,1]^{\s}$, where $a \in \Z^{\s}$.
For convenience, we include the empty set in the counting
and assign value $1$ to the corresponding summand.

Let $D$ be an integer cell. To compare $N(K,D)$ to $\Sigma(K)$ on
a simple example, take $K$ to be an integer box, i.e. the product
of $n$ intervals with integer endpoints and lengths $a_i \ge 0$,
$i = 1, \ldots, n$. Then $N(K,D) = \prod_1^n \max (a_i, 1)$ and
$\Sigma(K) = \prod_1^n (a_i + 1)$. Thus
$$
2^{-n} \Sigma(K)  \le N(K,D)  \le  \Sigma(K).
$$
The lower bound being trivially true for
any convex body $K$, an upper bound of this type is in general
difficult to prove. This motivates the following conjecture.

\begin{conjecture}[Covering Conjecture]
  Let $K$ be a convex body in $\R^n$ and $D$ be an integer cell.
  Then
  \begin{equation}                               \label{conj estimate}
  N(K, D) \le  \Sigma(C K)^C.
  \end{equation}
\end{conjecture}

Our main result is that the Covering Conjecture holds for a
body $D$ slightly larger that an integer cell, namely for
\begin{equation}                                 \label{tower D}
D = \Big\{ x \in \R^n \;:\;
           \frac{1}{n} \sum_1^n \exp \exp |x(i)| \le 3
    \Big\}.
\end{equation}
Note that the body $5D$ contains an integer cell
and the body $(5 \log \log n)^{-1} D$ is contained in an integer cell.

\begin{theorem}                                    \label{covering theorem}
  Let $K$ be a convex body in $\R^n$ and $D$ be the body \eqref{tower D}.
  Then
  $$
  N(K, D) \le  \Sigma(C K)^C.
  $$
\end{theorem}

As a useful consequence, the Covering Conjecture holds
for $D$ being an ellipsoid. This will follow by a standard
factorization technique for the absolutely summing operators.

\begin{corollary}                                 \label{ellipsoid}
  Let $K$ be a convex body in $\R^n$ and $D$ be
  an ellipsoid in $\R^n$ that contains an integer cell.
  Then
  $$
  N(K, D) \le  \Sigma(C K)^2.
  $$
\end{corollary}

As for the Covering Conjecture itself, it holds under the
assumption that the covering number is exponentially large in $n$.
Say, assume $N(K, D) \ge \exp(a n)$, where $D$ is an integer cell.
Then
\begin{equation}                                       \label{cube}
N(K,D)  \le  \Sigma(CK)^M, \ \ \ \text{where} \ \ M \sim
\log^{0.001} (2/a).
\end{equation}
The exponent $0.001$ can be replaced by arbitrarily small
positive number. This result also follows from Theorem \ref{covering
theorem}.

The usefulness of Theorem \ref{covering theorem} is understood
through a relation between the cell content and the combinatorial
dimension. Let $F$ be a class of real valued functions on a finite
set $\Omega$, which we identify with $\{1,\ldots, n\}$. Then we
can look at $F$ as a subset of $\R^n$ via the map $f \mapsto
(f(i))_{i=1}^n$. For simplicity assume that $F$ is a convex set;
the general case will not be much more difficult. It is then easy
to check that the combinatorial dimension $v := \vc(F,1)$ equals
exactly the maximal rank of a coordinate projection $P$ in $\R^n$
such that $PF$ contains a translate of the unit cube $P[0,1]^n$.
Then in the sum \eqref{lattice content} for the lattice content
$\Sigma(F)$, the summands with $\rank P > v$ vanish. The number of
nonzero summands is then at most $\sum_{k=0}^v \binom{n}{k}$.
Every summand is clearly bounded by $\vol(PF)$, a quantity which
can be easily estimated if the class $F$ is a priori
well bounded. So $\Sigma(F)$ is essentially bounded by
$\sum_{k=0}^v \binom{n}{k}$, and is thus controlled by the
combinatorial dimension $v$. This way, Theorem \ref{covering
theorem} or one of its consequences can be used to bound the
entropy of $F$ by its combinatorial dimension. Say, \eqref{cube}
implies \eqref{i:by linfty} in this way.

In some cases, $n$ can be removed from the bound on the entropy,
thus giving an estimate independent of the size of the domain
$\Omega$. Arguably the most general situation when this happens is
when $F$ is bounded in some norm and the entropy is computed with
respect to a weaker norm.
The entropy of the class $F$ with respect to a norm of a general
function space $X$ on $\Omega$ is
\begin{equation}                    \label{D(F,X,t)}
  D(F,X,t) = \log \sup \Big( n \mid \exists f_1, \ldots, f_n \in F \ \
  \forall i < j \  \|f_i-f_j\|_X \ge t \Big).
\end{equation}
Koltchinskii-Pollard entropy is then $D(F,t) = \sup_\mu D(F,
L_2(\mu), t)$, where the supremum is over all probability
measures supported by finite sets.
With the geometric representation as above,
\begin{equation}            \label{entropy vs packing}
  D(F,X,t) = \log N_{\rm pack} \Big( F, \frac{t}{2} \Ball(X) \Big)
\end{equation}
where $\Ball(X)$ denotes the unit ball of $X$
and $N_{\rm pack}(A,B)$ is the packing number,
which is the maximal number of disjoint
translates of a set $B \subseteq \R^n$ by vectors from a
set $A \subseteq \R^n$. The packing and the covering numbers
are easily seen to be equivalent,
\begin{equation}                                \label{covering vs packing}
  N_{\rm pack} (A, B)
  \le  N (A, B)
  \le  N_{\rm pack} (A, \frac{1}{2} B).
\end{equation}

To estimate $D(F, X, t)$, we have to be able to quantitatively
compare the norms in the function space $X$ an in another function
space $Y$ where $F$ is known to be bounded.
We shall consider Lorentz spaces, for which such a comparison
is especially transparent.
The {\em Lorentz space} $\L_\phi = \L_\phi(\Omega, \mu)$
is determined by its generating function $\phi(t)$,
which is a real convex function on $[0,\infty)$,
with $\phi(0)=0$, and increasing to infinity.
Then $\L_\phi$ is the space of functions $f$ on $\Omega$ such that
there exists a $\l > 0$ for which
\begin{equation}                                    \label{lorentz norm}
\mu \{ |f/\l| \ge t\}  \le  \frac{1}{\phi(t)}
\ \ \ \text{for all $t > 0$.}
\end{equation}
The norm of $f$ in $\L_\phi$ is the infimum of $\l > 0$ satisfying
\eqref{lorentz norm}.
Given two Lorentz spaces $\L_\phi$ and $\L_\psi$,
we look at their {\em comparison function}
$$
(\phi | \psi)(t) = \sup \{ \phi(s) \ | \  \phi(s) \ge \psi(ts) \}.
$$

Under the normalization assumption $\phi(1) = \psi(1) = 1$
and a mild regularity assumption on $\phi$ we prove
the following.
If a class $F$ is $1$-bounded in $\L_\psi$ then for all 
$0 < t < 1/2$
\begin{equation}                    \label{i:by lorentz}
D(\F, \L_\phi, t)  \le  C \; \vc(F, ct) \cdot \log (\phi|\psi)(t/2).
\end{equation}
An important point here is that the entropy is independent of
the size of the domain $\Omega$.
To prove \eqref{i:by lorentz}, we first perform a probabilistic
selection, which reduces the size of $\Omega$,
and then apply Theorem \ref{covering theorem}, in
which we replace $D$ by a larger set $\Ball(\L_\phi)$.

Of particular interest are the generating functions $\phi(t) = t^p$
and $\psi(t) = t^q$ with $1 \le p < q \le \infty$.
They define the weak $L_p$ and $L_q$ spaces respectively.
Their comparison function is $(\phi|\psi)(t) = t^{pq/(p-q)}$.
Then passing to usual $L_p$ spaces (which is not difficult)
one obtains from \eqref{i:by lorentz} the following.
If $F$ is $1$-bounded in $L_q(\mu)$ then for all $0 < t < 1/2$
\begin{equation}                                      \label{Lp}
D(F, L_p(\mu), t)
\le  C_{p,q} \; \vc(F,c_{p,q}t) \cdot \log (1/t),
\end{equation}
where $C_{p,q}$ and $c_{p,q} > 0$ depend only on $p$ and $q$.

First estimates of type \eqref{Lp} go back to the influental
works of Vapnik and Chervonenkis. In the main combinatorial
lemma of \cite{VC 81}, the volume of uniformly bounded convex
class was estimated via a quantity somewhat weaker than the
combinatorial dimension. Since we always have
$N(K,D) \ge \vol(K)/\vol(D)$,
the Vapnik-Chervonenkis bound is an asymptotically weaker
form of \eqref{Lp} for $p=2$ (say) and $q=\infty$.
Talagrand \cite{T 87, T 02} proved \eqref{Lp} for $p=2, q=\infty$
up to a factor of $\log^M (1/t)$ in the right side
and under minimal regularity (essentially under \eqref{i:min reg}).
Based on the method of N.Alon et al. from \cite{ABCH},
Bartlett and Long \cite{BL} proved \eqref{Lp} for $p=1, q=\infty$
with an additional factor of $\log(|\Omega|/vt)$ in the right
side, where $v = \vc(F, ct)$. The ratio $|\Omega|/v$ was removed from
this factor by Bartlett, Kulkarni and Posner \cite{BKP}, thus
yielding \eqref{Lp} with $\log^2(1/t)$ for $p=1, q=\infty$.
The optimal estimate \eqref{Lp} for all $p$ and for $q=\infty$
was proved by Mendelson and the second author as
the main result of \cite{MV}. FInally, the present paper
proves \eqref{Lp} for all $p$ and $q$. 

Finally, Theorems \ref{i:integral=integral} and
\ref{i:entropy=dimension} are proved by iterating \eqref{Lp}
with $2p = q \to \infty$ to get rid of both the logarithmic factor
and any boundedness assumptions.

\section{Covering by the Tower}                \label{s:by tower}

Fix a probability space $(\Omega, \mu)$. As most of our problems
have a discrete nature, they essentially reduce by
approximation to $\Omega$ finite and $\mu$ the uniform measure.
The core difficulties arise already in this finite setting,
although it took some time to fully realize this
(see \cite{T 96}).
This way we shall totally ignore measurability issues.

\paragraph{Tower}
Our main covering result works for a body in $\R^n$ which is $\log
\log n$ apart from the unit cube, while for the cube itself it
remains an open problem. This body is the unit ball of the Lorentz
space with generating function of the order $e^{e^t}$. For an extra
flexibility, we shall allow a parameter $\a \ge 2$, generally a
large number. The Lorentz space generated by the function
$$
\theta(t) = \theta_\a (t) = e^{\a^t - \a}, \ \ \ \ t \ge 1
$$
is called space the {\em tower space} and its unit ball
is called the tower. Since $\theta(1)=1$, it does not matter
how we define $\theta(t)$ for $0<t<1$ as long as
$\theta(0)=0$ and $\theta$ is convex; say, $\theta(t)=t$ will work.

In the discrete setting, we look at $\Omega$ being $\{1,
\ldots, n\}$ with the uniform probability measure $\mu$ on $\Omega$.
The tower space can be realized on $\R^n$ by identifying a
function on $\Omega$ with a point in $\R^n$ via the map
$f \mapsto (f_i)_{i=1}^n$.
The tower is then a convex symmetric body in $\R^n$,
and we denote it by $\tower^\a$. This body is equivalently
described by \eqref{tower D},
\[
c_1(\a) D  \subseteq  \tower^\a  \subseteq c_2(\a) D
\]
where positive $c_1(\a)$ and $c_2(\a)$ depend only on $\a$.

\paragraph{Coordinate convexity}
We stated our results for convex bodies but not necessarily
convex function classes. Convexity indeed plays very little role
in our work and is replaced by a much weaker notion of {\em
coordinate convexity}. This notion was originally motivated by
problems of calculus of variations, partial differential equations
and probability. The interested reader may consult the paper
\cite{M} and the bibliography cited there as an introduction to
the subject.

One can obtain a general convex body in $\R^n$ by cutting off
half-spaces.
Similarly, a general coordinate convex body in $\R^n$ is obtained
by cutting off octants, that is translates of the subsets of
$\R^n$ consisting of points with fixed and nonzero signs of the
coordinates. The {\em coordinate convex hull} of a set $K$ in
$\R^n$, denoted by $\cconv(K)$, is the minimal coordinate convex
set containing $K$. In other words, $\cconv(K)$ is what remains in
$\R^n$ after removal all octants disjoint from $K$. Clearly, every
convex set is coordinate convex; the converse is not true, as
shows the example of a cross $\{ (x,y) \; | \; x = 0 \text{ or } y
= 0 \}$ in $\R^2$.

\begin{center} 
\begin{picture}(100,100)(0,0)
\put(20,10) {\line(0,1){20}} \put(30,30) {\line(0,1){20}}
\put(0,30)  {\line(1,0){30}} \put(30,50) {\line(1,0){20}}
\put(50,50) {\line(0,1){20}} \put(50,70) {\line(1,0){40}}
\put(70,70) {\line(0,1){20}}
\end{picture} \\ \nopagebreak
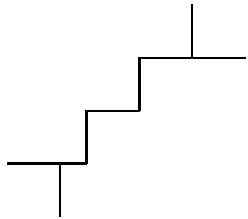
Example of a coordinate convex body in $\R^2$
\end{center}

\paragraph{Covering by the tower}
Let $A$ be a nonempty set in $\R^n$. In contrast to what
happens in classical convexity, a coordinate projection of a
coordinate convex set is not necessarily coordinate convex (a
pair of generic points in the plane is an example). Define the
{\em cell content} of $A$ as
$$
\Sigma(A) = \sum_P \text{number of integer cells in $\cconv(PA)$}
$$
where the sum is over all $2^n$ coordinate projections in $\R^n$,
including one $0$-dimensional projection,
for which the summand is set to be $1$. In many
applications $A$ will be a convex body, in which case $\cconv(PA)
= PA$.  The following is the main result of this section.

\begin{theorem}                                   \label{by tower}
  For every set $F$ in $\R^n$ and $\a \ge 2$,
  $$
  N(F, \tower^\a)  \le  \Sigma(C F)^\a
  $$
  where $C$ is an absolute constant.
\end{theorem}
It is plausible that the $\tower^\a$ can be replaced by the unit cube,
with $\a$ replaced by an absolute constant in the right hand side;
this is a slightly stronger version of the Covering Conjecture for
coordinate convex sets.

The proof of Theorem \ref{by tower}, which is a development upon
\cite{MV}, occupies next few subsections.

\paragraph{Separation on one coordinate}
Fix a set $F$ in $\R^n$ which contains more than one point.
Using \eqref{covering vs packing}, we can
find a finite subset $A' \subset F$ of cardinality
$N(F,\tower^\a)$ such that
no pair of points from $A'$ lies in a common translate of
$\frac{1}{2} \tower^\a$. Denote $A=2A'$. Then
$$
\forall x, y \in A, \; x \ne y: \ \ \
\|x - y\|_{\tower^\a}  \ge  1.
$$
Thus for a fixed pair $x \ne y$ there exists a $t > 0$ such that
$\mu \{ |x-y| > t \}  \ge  \frac{1}{\theta(t)}$.
Since $\theta(t) < 1$ for $t < 1$, we necessarily have $t \ge 1$, hence
$$
\exists t > 0: \ \ \ \mu \{ |x-y| > t \}  \ge  \frac{1}{\theta_0(t)}
$$
where
$$
\theta_0(t) = e^{\a^t - \a}, \ \ \ \ t \ge 0.
$$
By Chebychev's inequality,
$$
\E_i \; \theta_0 ( |x(i) - y(i)| )  \ge  1,
$$
where $\E_i$ is the expectation according to the uniform
distribution of the coordinate $i$ in $\{1,\ldots, n\}$. Let $x$
and $y$ be random points drawn from $A$ independently and
according to the uniform distribution on $A$. Then $x \ne y$ with
probability $1 - |A|^{-1} \ge \frac{1}{2}$, and taking the
expectation with respect to $x$ and $y$, we obtain
$$
\E_{x,y}\, \E_i \; \theta_0 ( |x(i) - y(i)| )  \ge  \frac{1}{2}.
$$
Changing the order of the expectation, we find a realization of
the random coordinate $i$ for which
\begin{equation}            \label{xi yi}
\E_{x,y}\; \theta_0 ( |x(i) - y(i)| )  \ge  \frac{1}{2}.
\end{equation}
Fix this realization.

Recall that a median of a real valued random variable $\xi$
is a number $M$ satisfying $\P (\xi \le M) \ge
1/2$ and $\P (\xi \ge M) \ge 1/2$. Unlike the expectation, the
median may be not uniquely defined.
We can replace $y(i)$ in \eqref{xi yi}
by a median of $x(i)$ using the following standard observation.

\begin{lemma}
  Let $\phi$ be a convex and nondecreasing function on $[0,\infty)$.
  Let $X$ and $Y$ be identically distributed random variables. Then
  $$
  \inf_a \E\, \phi(|X - a|)
  \le  \E\, \phi(|X - Y|)
  \le \inf_a \E\, \phi(2|X - a|).
  $$
\end{lemma}

\proof The first inequality follows from Jensen's inequality with
$a = \E X = \E Y$. For the second one, the assumptions on $\phi$
imply through the triangle and Jensen's inequalities that for
every $a$
$$
\phi(|X - Y|)  \le  \phi(|X - a| + |Y - a|) \le  \frac{1}{2}
\phi(2|X - a|) + \frac{1}{2} \phi(2|Y - a|).
$$
Taking the expectations on both sides completes the proof.
\endproof

Denote by $M$ a median of $x(i)$ over $x \in A$. We conclude
that
\begin{equation}                \label{expectation large}
\E_x\; \theta_0 ( 2 |x(i) - M| )  \ge  \frac{1}{2}.
\end{equation}

\begin{lemma}[Separation Lemma]
  Let $X$ be a random variable with median $M$.
  Assume that for every real $a$
  $$
  \P \{ X \le a \}^{1/\a} + \P \{ X > a+1 \}^{1/\a}  \le  1.
  $$
  Then
  $$
  \E \; \theta_0 (c |X - M| )  <  \frac{1}{2}.
  $$
\end{lemma}
In particular, the conclusion implies that the tower norm of the
random variable $X - M$ is bounded by an absolute constant.

\qquad

\proof One can assume that $M = 0$. With the notation $p(a) = \P
\{ X > a \}$, the assumption of the lemma implies that for every
$a$
$$
(1 - p(a)) + (p(a+1))^{1/\a} \le (1 - p(a))^{1/\a} +
(p(a+1))^{1/\a}  \le  1,
$$
hence
$$
p(a+1)  \le  p(a)^{1/\a}, \ \ \ a \in \R.
$$
Applying this estimate successively and using $p(0) =1-\P(x \le 0)
\le \frac{1}{2}$, we obtain $p(k)  \le  2^{-\a^k},  \ \ \ k \in
\N$. Then for every real number $a \ge 2$
$$
p(a)  \le  p([a])  \le  2^{-\a^{[a]}} \le  2^{-\a^{a-1}} \le
2^{-\a^{a/2}}.
$$
Repeating this argument for $-X$, we conclude that
$$
\P \{ |X| > a \}  \le  2^{1-\a^{a/2}}, \ \ \ a \ge 2.
$$
Then
$$
\P \{ e^{\a^{c|X|}} > s \}  \le  2^{1 - (\log s)^{1/2c}} \le
2s^{-\a^{1-2c}}, \ \ \ s \ge e^{\a^{2c}}.
$$
Integrating by parts and using this tail estimate, we have
\begin{align*}
\E\; \theta_0(c|X|) &=   e^{-\a} \E e^{\a^{c|X|}} \le  e^{-\a} \Big[
e^{\a^{2c}} + \int_{e^{\a^{2c}}}^\infty
                   2s^{-\a^{1-2c}} \; ds \Big] \\
&=   e^{-\a + \a^{2c}}  +  2 (\a^{1-2c} - 1)^{-1} e^{-2\a +
\a^{2c}} =: h(\a, c).
\end{align*}
For a fixed $c \le 1/4$, the function $h(\a, c)$ decreases as a
function of $\a$ on $[2,\infty)$, and $h(2, 0) = e^{-1} + 2e^{-3}
\approx 0.47 < \frac{1}{2}$. Hence for a suitable choice of the
absolute constant $c > 0$,
$$
h(\a, c)  \le  h(2, c)  < \frac{1}{2}
$$
because $\a \ge 2$. This completes the proof.
\endproof

\qquad

Applying the Separation Lemma to the random variable $\frac{2}{c}
x(i)$ together with \eqref{expectation large},
we find an $a \in \R$ so that
$$
\mu \{ x(i) \le a \}^{1/\a} + \mu \{ x(i) > a + c \}^{1/\a} > 1,
$$
where $\mu$ is the uniform measure on $A$. Equivalently, for the
subsets $A_-$ and $A_+$ of $A$ defined as
\begin{equation}                \label{sons}
A_- = \{ x \;:\; x(i) \le a \},  \ \ \ A_+ = \{ x \;:\; x(i) >  a
+ c \}
\end{equation}
we have
\begin{equation}                \label{sons large}
|A_-|^{1/\a} + |A_+|^{1/\a}  >  |A|^{1/\a}.
\end{equation}
Here $|A|$ denotes the cardinality of the set $A$.

\paragraph{Separating tree}
This and the next step are versions of corresponding steps of
\cite{MV}, where they were written in terms of function classes.
Continuing the process of separation for each $A_-$ and $A_+$, we
construct a {\em separating tree} of subsets of $A$.

A tree of nonempty subsets of a set $A$ is a finite collection $T$ of
subsets of $A$ such that every two elements in $T$ are either
disjoint or one contains the other. A {\em son} of an element $B
\in T$ is a maximal (with respect to inclusion) proper subset of
$B$ which belongs to $T$. An element of with no sons is called a
{\em leaf}, an element which is not a son of any other element is
called a {\em root}.

\begin{definition}
  Let $A$ be a class of functions on $\Omega$ and $t > 0$.
  A {\em $t$-separating tree $T$ of $A$} is a tree of subsets of $A$
  whose only root is $A$ and such that every element $B \in T$
  which is not a leaf has exactly two sons $B_+$ and $B_{-}$
  and, for some coordinate $i \in \Omega$,
  $$
  f(i) \ge g(i) + t
  \ \ \ \text{for all $f \in B_+$,  $g \in B_-$.}
  $$
\end{definition}
If $|A_-| > 1$, we can repeat the separation on one coordinate for
$A_-$ (note that this coordinate may be different from $i$). The
same applies to $A_+$. Continuing this process of separation until
all the resulting sets are singletons, we arrive at

\begin{lemma}               \label{large tree}
  Let $A \subset \R^n$ be a finite set
  whose points are $1$-separated in $\tower^{\a}$-norm.
  Then there exists a $c$-separating tree of $A$ with at least
  $|A|^{1/\a}$ leaves.
\end{lemma}
This separating tree improves in a sense the set $A$ which was
already separated. Of course, the leaves in this tree are
$c$-separated in the $L_\infty$-norm, but the tree also shows some
pattern in the coordinates on which they are separated. This will
be used in the next section where we further improve the
separation of $A$ by constructing in it many copies of a discrete
cube (on different subsets of coordinates).

However note that the assumption on $A$, that it is separated in
the tower norm, is stronger than being separated in the
$L_\infty$-norm.

\qquad

\proof
We proceed by induction on the cardinality of $A$.
The claim is trivially true for singletons.
Assume that $|A| > 1$ and that the claim holds for all
sets of cardinality smaller than $|A|$.
By the separation procedure described above, we can
find two subsets $A_-$ and $A_+$ satisfying \eqref{sons}
and \eqref{sons large}.
The strict inequality in \eqref{sons large} implies that
the cardinalities of both sets is strictly smaller than $|A|$.
By the induction hypothesis, both $A_-$ and $A_+$ have
$c$-separating trees $T_-$ and $T_+$ with at least
$|A_-|^{1/\a}$ and $|A_+|^{1/\a}$ leaves respectively.

Now glue the trees $T_-$ and $T_+$ into one tree $T$
of subsets of $A$ by declaring $A$ the root of $T$
and $A_-$ and $A_+$ the sons of $A$. By \eqref{sons},
$f(i) \ge g(i) + c$ for all $f \in A_+$, $g \in A_-$.
Therefore $T$ is a $c$-separating tree of $A$.
The number of leaves in $T$ is the sum of the
number of leaves of $T_-$ and $T_+$, which is at least
$|A_-|^{1/\a} + |A_+|^{1/\a}  >  |A|^{1/\a}$ by
\eqref{sons large}.
This proves the lemma.
\endproof

\paragraph{Coordinate convexity and counting cells}
Recall that $|A|=N(F, \tower^{\a})$.
 We shall prove the following fact which, together with Lemma
\ref{large tree}, finishes the proof.

\begin{lemma}                              \label{many cells}
  Let $A$ be a set in $\R^n$, and $T$ be a $2$-separating tree of
  $A$. Then
  $$
  \text{Number of leaves in $T$}  \le  \Sigma(A).
  $$
\end{lemma}

The value $2$ is exact here. For example, the open cube
$A = (-1,1)^n$ has $\Sigma(A) = 1$, because $A$ contains
no integer cells. However, for every $\e > 0$
one easily constructs a $(2-\e)$-separating tree of $A$
with $2^n$ leaves.

We ask what it means for a cell to be contained in the
coordinate convex hull of a set.  A cell $\CC$ in $\R^n$ defines
$2^n$ octants in a natural way. Let $\theta \in \{-1,1\}^n$ be a
choice of signs. A closed octant with the vertex $z \in \R^n$ is
the set
$$
\mathcal{O}_{\theta}(z) = \{ x=(x_1, \ldots x_n) \in \R^n \mid
(x_i-z_i) \cdot \theta_i \ge 0 \quad \mbox{for } i=1, \ldots n\}.
$$
 The octants generated by a cell
are those who have only one common point with it (a vertex).

\begin{lemma}                                                 \label{cell}
  Let $A$ be a set in $\R^n$ and $\CC$ be a cell of $\Z^n$.
  Then $\CC \subset \cconv(A)$ if and only if
  $A$ intersects all the octants generated by $\CC$.
\end{lemma}
The proof is straightforward and we omit it.
\endproof

\qquad

\noindent {\bf Proof of Lemma \ref{many cells}.}
It will suffice to prove that
\begin{equation}            \label{three sigmas}
  \text{if $A_-$ and $A_+$ are the sons of $A$, then} \ \ \Sigma(A_-) +
  \Sigma(A_+)  \le  \Sigma(A).
\end{equation}
Indeed, assuming that \eqref{three sigmas} one can complete
the proof by induction on the cardinality of $A$ as follows.
The lemma is trivially true for singletons.
Assume that $|A| > 1$ and that the lemma holds for all
sets of cardinality smaller than $|A|$.
Let $A_-$ and $A_+$ be the sons of $A$.
Define $T_-$ to be the colection of sets from $T$
that belong to $A_-$; then $T_-$ is a separating tree of $A_-$.
Do similarly for $T_+$.
Since both $A_-$ and $A_+$ have cardinalities smaller than $|A|$,
the induction hypothesis applies to them.
Hence by \eqref{three sigmas} we have
\begin{align*}
\Sigma(A) \ge \Sigma(A_-) + \Sigma(A_+)
 &\ge  \text{(number of leaves in $T_-$)}
    + \text{(number of leaves in $T_+$)} \\
 &= \text{number of leaves in $T$}.
\end{align*}
This proves the lemma, so the only remaining thing is to prove
\eqref{three sigmas}.

In the proof of \eqref{three sigmas}, when it creates no confusion,
we will denote by $\Sigma(A)$ not only the cardinality,
but also the set of all pairs $(P, \CC)$ for which
$\CC \subset \cconv(PA)$. For this to
be consistent, we introduce a $0$-dimensional cell $\emptyset$,
and always assume that the $0$-dimensional projection
along with the empty cell are in $\Sigma(A)$ provided $A$ is nonempty.

Clearly, $\Sigma(A_-) \cup \Sigma(A_+)  \subseteq  \Sigma(A)$. To
complete the proof, it will be enough to construct an injective
mapping $\Phi$ from $\Sigma(A_-) \cap \Sigma(A_+)$ into $\Sigma(A)
\setminus (\Sigma(A_-) \cup \Sigma(A_+))$. We will do this by
gluing identical cells from $\Sigma(A_-) \cap \Sigma(A_-)$ into a
larger cell; this idea goes back to \cite{ABCH}.

Fix a pair $(P, \CC) \in \Sigma(A_-) \cap \Sigma(A_+)$. Without
loss of generality, we may assume that $A_-$ and $A_+$ are
$2$-separated on the first coordinate. Then there exists an
integer $a$ such that
\begin{equation}                                 \label{1-separated}
  \text{$x(1) \le a$   for $x \in A_-$,} \ \ \
  \text{$x(1) \ge a+1$ for $x \in A_+$.}
\end{equation}
The coordinate projection $P$ must annihilate the first
coordinate, otherwise \eqref{1-separated} would imply that the
sets $PA_-$ and $PA_+$ are disjoint, which would contradict to our
assumption that their coordinate convex hulls both contain the
cell $\CC$.

{\bf Trivial case: $\rank P = 0$.} In this case, let $P'$ be the
coordinate projection that annihilates all the coordinates except
the first. Since both $A_-$ and $A_+$ are nonempty, $P'A$ contains
points for which $x(1) \le a$ and $x(1) \ge a+1$. Hence
$\cconv(P'A)$ contains the one-dimensional cell $\CC' = [a, a+1]$.
So, we can define the action of $\Phi$ on the trivial pair as
$\Phi : (P, \emptyset) \mapsto (P', \CC')$.

{\bf Nontrivial case: $\rank P > 0$.} Without loss of generality
we may assume that $P$ retains the coordinates $\{2, 3, \ldots,
k\}$ with some $2 \le k \le n$, and annihilates the others. Let
$P'$ be the coordinate projection onto $\R^k$, so $\CC' = [a, a+1]
\times \CC$ is a cell in $\R^k$. We claim that $(P', \CC') \in
\Sigma(A)$. By the assumption, the cell $\CC$ lies in both
$\cconv(PA_-)$ and $\cconv(PA_+)$ . In light of Lemma \ref{cell},
$PA_-$ and $PA_+$ each intersect all the octants generated by
$\CC$, and we need to show that $PA'$ intersects any octant $\OO'$
generated by $\CC'$. This octant must be of the form either $\OO'
= \{ x \in \R^k \;:\; x(1) \le a, Px \in \OO \}$ or $\OO' = \{ x
\in \R^k \;:\; x(1) \ge a+1, Px \in \OO \}$, where $\OO$ is some
octant generated by the cell $\CC$. Assume the second option
holds. Pick a point $z \in A_+$ such that $Pz \in PA_+ \cap \OO$.
Then $P'z(1) = z(1) \ge a+1$, so $P'z \in P'A_+ \cap \OO'$. A
similar argument (with $A_-$) works if $\OO$ is of the first form.
This proves the claim, and we again define the action of $\Phi$ as
$\Phi : (P, \CC) \mapsto (P', \CC')$.

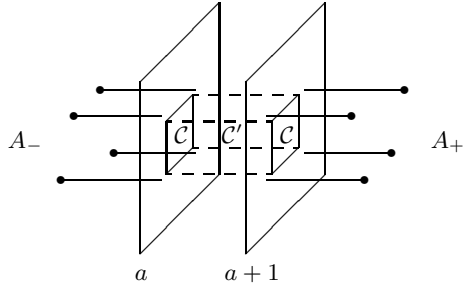
\begin{center} 
\begin{picture}(120,120)(0,-20)
\multiput(30,0)(40,0){2} {\line(1,1){30}} \multiput(30,0)(40,0){2}
{\line(0,1){65}} \multiput(30,65)(40,0){2} {\line(1,1){30}}
\multiput(60,30)(40,0){2} {\line(0,1){65}}
\multiput(40,30)(40,0){2} {\line(1,1){10}}
\multiput(40,30)(40,0){2} {\line(0,1){20}}
\multiput(40,50)(40,0){2} {\line(1,1){10}}
\multiput(50,40)(40,0){2} {\line(0,1){20}}
\multiput(40,30)(10,10){2} {\dashbox{4}(40,20)}
\put(0,28) {\circle*{3}}
  \put(0,28) {\line(1,0){38}}
\put(5,52) {\circle*{3}}
  \put(5,52) {\line(1,0){33}}
\put(15,62) {\circle*{3}}
  \put(15,62) {\line(1,0){36}}
\put(20,38) {\circle*{3}}
  \put(20,38) {\line(1,0){31}}
\put(115,28) {\circle*{3}}
  \put(78,28) {\line(1,0){37}}
\put(110,52) {\circle*{3}}
  \put(78,52) {\line(1,0){32}}
\put(130,62) {\circle*{3}}
  \put(92,62) {\line(1,0){38}}
\put(125,38) {\circle*{3}}
  \put(92,38) {\line(1,0){33}}
\multiput(43,42)(40,0){2}{\footnotesize{$\CC$}}
\put(61,42){\footnotesize{$\CC'$}}
\put(28,-10){\footnotesize{$a$}}
\put(62,-10){\footnotesize{$a+1$}}
\put(-20,40){\footnotesize{$A_-$}}
\put(140,40){\footnotesize{$A_+$}}
\end{picture}

Nontrivial case: Gluing two copies of $\CC$ into a larger cell
$\CC'$
\end{center}

To check that the range of $\Phi$ is disjoint from both
$\Sigma(A_-)$ and $\Sigma(A_+)$, assume that the pair $(P', \CC')$
constructed above is in $\Sigma(A_-)$. This means that $\CC'$ lies
in $\cconv(QA_-)$ for some coordinate projection $Q$. This
projection must retain the first coordinate because the cell
$\CC'$ is non-degenerating on the first coordinate by its
construction. Therefore, since $x(1) \le a$ for all $x \in A_-$,
the same must hold for all $x \in Q(A_-)$, and hence also for all
$x \in \cconv(QA_-)$. On the other hand, there clearly exist
points in $\CC'$ with $x(1) = a+1 > a$. Hence $\CC'$ can not lie
in $\cconv(QA_-)$. A similar argument works for $A_+$. Therefore
the range of $\Phi$ is as claimed.

Finally, $\Phi$ is trivially injective because the map $\CC
\mapsto \CC'$ is injective.
\endproof

\qquad

Theorem \ref{by tower} follows from Lemma \ref{large tree} and
Lemma \ref{many cells}.

\qquad

\remark The proof does not use the fact that the probability
measure on $\Omega = \{1, \ldots, n\}$, underlying the tower
space, is uniform. In fact, Theorem \ref{by tower} holds for any
probability measure on $\{1, \ldots, n\}$. This will help us in
next section.

\section{Covering by Ellipsoids and Cubes}
                 \label{s:by ellipsoids and cubes}

The Covering Conjecture holds if we cover by ellipsoids containing
the unit cube rather by the unit cube itself. This nontrivial fact
is a consequence of Theorem \ref{by tower}.

\begin{theorem}                                       \label{by ellipsoids}
  Let $A$ be a set in $\R^n$ and $D$ be an ellipsoid containing
  the cube $[0,1]^n$. Then
  $$
  N(A, D)  \le  \Sigma(C A)^2
  $$
  where $C$ is an absolute constant.
\end{theorem}

This result will be used in Section \ref{s:sections} to find nice
sections of convex bodies.

\qquad

\proof Translating the ellipsoid $D$, we can assume that $2D$
contains the cube $[-1,1]^n$, which is the unit ball of the space
$l_\infty^n$. Call $X$ the normed space $(\R^n,
\left\|\cdot\right\|_{2D})$. Then $X$ is isometric to $l_2^n$. Let
$T : l_\infty^n \to X$ be the formal identity map and $S : X \to
l_2^n$ be an isometry. Finally, define $u = S T : l_\infty^n \to
l_2^n$ and note that $\|u\| \le 1$. Recall that every linear
operator $u : l_\infty^n \to l_2^n$ is $2$-summing and its
$2$-summing norm $\pi_2(u)$ satisfies $\pi_2(u)  \le  \sqrt{\pi/2}
\; \|u\|$, see \cite{TJ} Corollary 10.10. Thus $\pi_2(u)  \le
\sqrt{\pi/2}$. By Pietsch's factorization theorem (see \cite{TJ}
Theorem 9.3) there exists a probability measure $\mu$ on $\Omega =
\{1, \ldots, n\}$ such that for all $x \in \R^n$
$$
\|u x\|  \le  \sqrt{\pi/2} \; \|x\|_{L_2(\Omega,\mu)}.
$$
Since $\|u x\| = \|S^{-1} u x\|_X = \|T x\|_X = \|x\|_X$, we have
\begin{equation}                                    \label{X L2}
\frac{1}{\sqrt{\pi/2}} \|x\|_X
\le  \|x\|_{L_2(\Omega,\mu)}.
\end{equation}
On the other hand, the norm of the Lorentz space generated by
$\theta_2(t) = e^{2^{t} - 2}$ clearly dominates the $L_2$ norm:
for every $x \in \R^n$,
\begin{equation}                                    \label{L2 tower}
\|x\|_{L_2(\Omega, \mu)}  \le  C \|x\|_{\L_{\theta_2}(\Omega,\mu)}
\end{equation}
where $C$ is an absolute constant.
Denoting by $\tower^2(\mu)$ the unit ball of the norm
in the right hand side of \eqref{L2 tower}, we
conclude from \eqref{X L2} and \eqref{L2 tower} that
$$
\tower^2(\mu) \subseteq C' D
$$
where $C'$ is an absolute constant.
Then by Theorem \ref{by tower} and the remark after its proof,
$$
N(A, D)  \le  N(C'A, \tower^2(\mu))  \le  \Sigma(C'' A)^2
$$
where $C''$ is an absolute constant.
\endproof

Next theorem is a partial positive solution to the
Covering Conjecture itself. We prove the conjecture with a
mildly growing exponent.

\begin{theorem}                                        \label{by cubes}
  Let $A$ be a set in $\R^n$ and $\e > 0$.
  Then for the integer cell $Q = [0,1]^n$
  $$
  N(A, Q) \le   \Sigma(C \e^{-1} A)^M
  $$
  with $M = 4 \log^\e (e + n / \log N(A, Q))$,
  and where $C$ is an absolute constant.
\end{theorem}

In particular, this proves the Covering Conjecture in case when
the covering number is exponential in $n$: if $N(A, Q) \ge \exp(\l
n)$, $\l < 1/2$, then $M \le C \log^\e(1/\l)$.

\qquad

For the proof of the theorem, we first cover $A$ by towers, and
then towers by cubes. Formally,
\begin{align}                            \label{entropy product}
  N(A, Q)
  &\le  N(A, \e \tower^\a) \; N(\e \tower^\a, Q) \nonumber \\
  & =    N(\e^{-1} A, \tower^\a) \; N(\tower^\a, \e^{-1} Q).
\end{align}

\begin{lemma}                               \label{tower by cubes}
  For every $t \ge 4$,
  $$
  N(\tower^\a, tQ)  \le  \exp(C e^{-\frac{1}{4}\a^{t/2}} n)
  $$
  where $C$ is an absolute constant.
\end{lemma}

\proof We count the integer points in the tower. For $x \in \R^n$,
define a point $x' \in \Z^n$ by $x'(i) = \sign(x(i)) [x(i)]$.
Every point $x \in \tower^\a$ is covered by the cube $x' +
[-1,1]^n$, so
\begin{align*}
N = N(\tower^\a, tQ)
  &= N(2t^{-1} \tower^\a, 2Q)
  \le |\{ x' \in \Z^n\ | \  x \in 2t^{-1} \tower^\a \} |\\
  &\le | 2t^{-1} \tower^\a \cap \Z^n|.
\end{align*}
For every $x \in 2t^{-1} \tower^\a \cap \Z^n$,
$$
|\{ i : |x(i)| = j \}|  \le  e^{-\a^{t j / 2} + \a} n =: k_j, \ \
\ j \in \N.
$$
Let $J$ be the largest number $j$ such that $k_j \ge 1$. Then
$$
N \le \prod_{j=1}^J \binom{n}{k_j} 2^{k_j},
$$
as for every $j$ there are at most $\binom{n}{k_j}$ ways to choose
the the level set $\{ i : |x(i)| = j \}$, and at most $2^{k_j}$
ways to choose signs of $x(i)$.

Let $\b_j = k_j / n$. Since $\a \ge 2$ and $t \ge 2$, $\b_j < 1/4$.
Then
$\binom{n}{k_j}
\le \big( e / \b_j \big)^{\b_j n}
\le \exp(C \b_j^{1/2} n)$.
Hence
$$
N \le \exp \Big( C_1 \sum_{j=1}^J \b_j^{1/2} n \Big) \le \exp(C_2
\b_1^{1/2} n) \le  \exp(C_2 e^{-\frac{1}{4}\a^{t/2}} n).
$$
This completes the proof.
\endproof

\qquad

\noindent{\bf Proof of Theorem \ref{by cubes}. } We can assume
that $0 < \e < c$ where $c > 0$ is any absolute constant.
We estimate the second factor in \eqref{entropy product}
by Lemma \ref{tower by cubes}. With $\a = M/2$,
\begin{align*}
N(\tower^\a, \e^{-1} Q)
  &\le  \exp \Big[ C \Big( e + \frac{n}{\log N(A, Q)}
                   \Big)^{-2^{1/2\e}/4} n \Big]\\
&\le  \exp \Big[ C e^{-2^{1/2\e}/4 + 1}
             \Big( \frac{n}{\log N(A, Q)} \Big)^{-1} n \Big] \\
&\le  N(A, Q)^{1/2}.
\end{align*}
Then \eqref{entropy product} and Theorem \ref{by tower} imply that
$$
N(A, Q)  \le  N(\e^{-1} A, \tower^{M/2})^2  \le \Sigma(c \e^{-1}
A)^M.
$$
The proof is complete.
\endproof

\qquad

Theorem \ref{by cubes} applies to a combinatorial problem
studied by N.Alon et al. \cite{ABCH}.

\begin{theorem}                               \label{by Linfty}
  Let $F$ be a class of functions on an $n$-point set $\Omega$
  with the uniform probability measure $\mu$.
  Assume $F$ is $1$-bounded in $L_1(\Omega, \mu)$. Then
  for $0 < \e < 1$ and for $0 < t < 1/2$
  \begin{equation}                            \label{eq by Linfty}
  D_\infty(F, t)
  \le  C v \, \log(n/vt) \cdot \log^\e(2n/v)
  \end{equation}
  where $v = \vc(F, c \e t)$.
\end{theorem}

N.Alon et al. \cite{ABCH} proved under a somewhat stronger
assumption ($F$ is $1$-bounded in $L_\infty$) that
\begin{equation}                                  \label{abch}
D_\infty(F, t) \le  C v \, \log(n/vt) \cdot \log(n/t^2),
\ \ \ \text{where $v = \vc(F, c t)$}.
\end{equation}
Thus $D_\infty(F, t) = O(\log^2 n)$.
It was asked in \cite{ABCH} whether the exponent $2$
can be reduced to some constant between $1$ and $2$.
Theorem \ref{by Linfty} answers this in positive.
It remains open whether the exponent can be made $1$.
A partial case of Theorem \ref{by Linfty},
for $\e=2$ and for uniformly bounded classes,
was proved in \cite{MV 02}.

It is important that, unlike in \eqref{abch},
the size of the domain $n$ appears in \eqref{eq by Linfty}
always in the ratio $n/v$. Assume, for
example, that one knows {\em a priori} that the entropy is large:
for some constant $0 < a < 1/2$
$$
D_\infty(F, t) \ge a n.
$$
Then by \eqref{eq by Linfty} we have
$a n  \le  C v \, \log(n/vt) \cdot \log^\e (2n/v)$.
Dividing by $n$ and solving for $n/v$, we get
$$
n/v  \le
\frac{C}{a} \Big[ \log \Big(\frac{1}{t} \Big)
  \log^\e \Big(\frac{1}{a} \log \frac{1}{t} \Big)
  + \log^{1+\e} \Big(\frac{1}{a} \Big)  \Big]
$$
and putting this back into \eqref{eq by Linfty} we obtain
$$
D_\infty(F, t)  \le  C v \, \log \Big(\frac{1}{a t} \Big)
  \cdot \log^\e \Big(\frac{1}{a} \log \frac{1}{t} \Big).
$$
We see that $n$, the size of the domain $\Omega$,
disappeared from the entropy estimate.
Such domain-free bounds, to which we shall return in the next
section, are possible only because $n$ enters into
the entropy estimate \eqref{eq by Linfty} in a ratio $n/v$.

\qquad

To prove Theorem \ref{by Linfty}, we identify the $n$-point domain
$\Omega$ with $\{1, \ldots, n\}$ and
realize the class of functions $F$ as a subset of $\R^n$
via the map $f \mapsto (f(i))_{i=1}^n$.
The geometric meaning of the combinatorial dimension of $F$
is then the following.

\begin{definition}
  The {\em combinatorial dimension} $\vc(A)$ of a set $A$ in $\R^n$
  is the maximal rank of a coordinate projection $P$
  in $\R^n$ so that $\cconv(PA)$ contains an integer cell.
\end{definition}
This agrees with the classical Vapnik-Chernovenkis definition for
sets $A \subseteq \{0,1\}^n$, for which $\vc(A)$ is defined as the
maximal rank of a coordinate projection $P$ such that
$PA = P(\{0,1\}^n)$.

\begin{lemma}                                           \label{to Rn}
  $\vc(F, 1) = \vc(F)$, where $F$ is treated as a function class
  in the left hand side and as a subset of $\R^n$ in the right
  hand side.
\end{lemma}

\proof By the definition, $\vc(F, 1)$ is the maximal cardinality of
a subset $\s$ of $\{1, \ldots, n\}$ which is $1$-shattered by $F$.
Being $1$-shattered means that there exists a point $h \in \R^n$
such that for every partition $\s = \s_- \cup \s_+$
one can find a pont $f \in F$ with $f(i) \le h(i)$ if $i \in \s_-$
and $f(i) \ge h(i) + 1$ if $i \in \s_+$.
This means exactly that $P_\s F$
intersects each octant generated by the cell $\CC = h + [0,1]^\s$,
where $P_\s$ denotes the coordinate projection in $\R^n$ onto
$\R^\s$. By Lemma \ref{cell} this means that $\CC \subset
\cconv(PF)$. Hence $\vc(F, 1) = \vc(F)$.
\endproof

For a further use, we will prove Theorem \ref{by Linfty}
under a weaker assumption, namely that $F$ is $1$-bounded in
$L_p(\mu)$ for some $0 < p < \infty$.
When $F$ is realized as a set in $\R^n$,
this assumption means that $F$ is a subset of the unit ball
of $L_p^n$, which is
$$
\Ball(L_p^n)
= \Big\{ x \in \R^n \;:\; \sum_1^n |x(i)|^p  \le  n \Big\}.
$$
We will apply to $F$ the covering Theorem \ref{by cubes} and then
estimate $\Sigma(F)$ as follows.

\begin{lemma}                                    \label{A in L1}
  Let $A$ be a subset of $a \cdot \Ball(L_p^n)$
  for some $a \ge 1$ and $0 < p \le \infty$. Then
  $$
  \Sigma(A)  \le  \Big( \frac{C_1(p) a n}{v} \Big)^{C_2(p) v}
  $$
  where $v = \vc(A)$, $C_1(p) = C(1 + \frac{1}{\sqrt{p}})$
  and $C_2(p) = 1 + \frac{1}{p}$.
\end{lemma}

\proof We look at
$$
\Sigma(A) = \sum_P \text{number of integer cells in $\cconv(PA)$}
$$
and notice that by Lemma \ref{to Rn}, $\rank P \le \vc(A) = v$ for
all $P$ in this sum. Since the number of integer cells in a set is
always bounded by its volume,
$$
\Sigma(A)  \le  \sum_{\rank P \le v} \vol(\cconv(PA)) \le
\sum_{\rank P \le v} \vol \Big(P(a \cdot \Ball(L_p^n)) \Big)
$$
where the volumes are considered in the corresponding subspaces
$P(\R^n)$. By the symmetry of $L_p^n$, the summands with the
same $\rank P$ in the last sum are equal. Then the sum equals
\begin{equation}                \label{1+sum}
1 + \sum_{k=1}^v \binom{n}{k} a^k \,
                  \vol_k \Big( P_k(\Ball(L_p^n)) \Big)
\end{equation}
where $P_k$ denotes the coordinate projection in $\R^n$ onto $\R^k$.
Note that $P_k(\Ball(L_p^n)) = (n/k)^{1/p} \, \Ball(L_p^k)$
and recall that $\vol(\Ball(L_p^k)) \le C_1(p)^k$,
see \cite{Pi} (1.18).
Then the volumes in \eqref{1+sum} are bounded by
$(n/k)^{k/p} C_1(p)^k  \le  (C_1(p) n / k)^{C_2(p) k}$.
The binomial coefficients in \eqref{1+sum} are estimated
via Stirling's formula as $\binom{n}{k} \le (en/k)^k$.
Then \eqref{1+sum} is bounded by
$$
1 + \sum_{k=1}^v \Big( \frac{e n}{k} \Big)^k
  a^k \Big( \frac{n}{k} \Big)^{k/p} C_1(p)^k
\le  \Big( \frac{C \cdot C_1(p) a n}{v} \Big)^{C_2(p) v}.
$$
This completes the proof.
\endproof

\qquad

\noindent {\bf Proof of Theorem \ref{by Linfty}. } Viewing $F$ as
a set in $\R^n$, we notice from \eqref{entropy vs packing} and
\eqref{covering vs packing} that
\begin{equation}                                        \label{Dinfty vs N}
D_\infty(F,t)  \le  \log N(F, 2tQ)  \le  D_\infty(F, t/2)
\end{equation}
where $Q = [0,1]^n$.
Therefore it is enough to estimate $N = N(F, 2tQ)$.
We apply successively the covering Theorem \ref{by cubes}
and Lemma \ref{A in L1} with $p=1$:
\begin{equation}                    \label{N small}
  N = N \big( \frac{1}{2t} F, Q \big)
  \le \Sigma \Big( \frac{C}{\e t} F \Big)^M
  \le \Big( \frac{C n}{\e t v} \Big)^{C M v}
\end{equation}
where $v = \vc(\frac{c}{\e t} F) = \vc(F, \frac{\e t}{c})$
and $M = 4 \log^\e (e + n / \log N)$.
Define the number $a > 0$ by $N = \exp(a n)$.
Then $M = 4 \log^\e (e + \frac{1}{a})$
and taking logarithms in \eqref{N small}
we have $a n  \le  C M v \, \log (\frac{C n}{\e v t})$.
Dividing by $M n$, we obtain
$$
\frac{a}{\log^\e (e + \frac{1}{a})} \le \frac{C v}{n} \log
\Big(\frac{C n}{\e v t} \Big).
$$
This implies
$$
a  \le  \frac{C v}{n} \log \Big( \frac{C n}{\e v t} \Big)
         \log^\e \Big( \frac{C n}{v}
           \big/ \log \Big(\frac{C n}{\e v t} \Big) \Big)
    \le  \frac{C v}{n} \log \Big(\frac{C n}{\e v t} \Big)
         \log^\e \Big( \frac{C n}{v} \Big)
$$
and multiplying by $n$ we obtain
\begin{equation}                                  \label{log N small}
\log N \le  C v \, \log(C n / v \e t) \cdot \log^\e(C n/v).
\end{equation}
It remains to remove $\e$ from the denominator by a routine argument.

Consider the function
$$
\phi(\e) =  \log^\e(C n / v), \ \ \ \text{where $v = v(\e)$ as before}.
$$
As $\e$ decreases to zero, $v(\e)$ increases, thus $\phi(\e)$ decreases
to $1$. Define $\e_0$ so that $\phi(\e_0) = e$.

{\bf Case 1.} Assume that $\e \ge \e_0$. Then $\phi(\e) \ge e$,
thus $\e \ge 1/\log \log (C n / v)$, so $C n / v \e t \le  (C n /
v t)^2$. Using this in \eqref{log N small} we obtain
\begin{equation}                                \label{log N final}
\log N  \le  C v \, \log(C n / v t) \cdot \log^\e(C n/v).
\end{equation}

{\bf Case 2.} Let  $\e < \e_0$. Then $\phi(\e) \le e$, so by
\eqref{log N small},
\begin{equation}                \label{case 2 small}
\log N \le  C v(\e_0) \, \log(C n / v(\e_0) \e_0 t) \cdot e.
\end{equation}
As in case 1, we have
$C n / v(\e_0) \e_0 t \le  (C n / v(\e_0) t)^2$.
Using this in \eqref{case 2 small}, we obtain
$$
\log N \le  C' v(\e_0) \, \log(C n / v(\e_0) t)
\le C' v \log(C n / v t),
$$
because $v(\e_0) \le v(\e) = v$.
In particular, we have \eqref{log N final} also in this case.
In view of \eqref{Dinfty vs N}, this completes the proof.
\endproof

\section{Covering by balls of Lorentz spaces}       \label{s:by lorentz}

So far we imposed no assumptions on the set $A \subset \R^n$
which we covered.
If $A$ happens to be bounded in some norm $\|\cdot\|$, a
new phenomenon occurs. The covering numbers of $A$ by balls in any
norm slightly weaker than $\|\cdot\|$ become independent of the
dimension $n$; the parameter that essentially controls them is the
combinatorial dimension of $A$.

This phenomenon is best expressed in the functional setting for
Lorentz norms \eqref{lorentz norm}, because they are especially
easy to compare. Given two generating functions $\phi$ and $\psi$,
we look at their {\em comparison function}
$$
(\phi | \psi)(t) = \sup \{ \phi(s) \;|\; \phi(s) \ge \psi(ts) \}.
$$
Fix a probability space $(\Omega, \mu)$. The comparison function
helps us measure to what extent the norm in $\L_\phi =
\L_\phi(\Omega, \mu)$ is weaker than the norm in $\L_\psi =
\L_\psi(\Omega, \mu)$.

Just for the normalization, we assume that
\begin{equation}                                      \label{normalization}
\phi(1) = \psi(1) = 1.
\end{equation}
Let $2 \le \a < \infty$. We rule out the extremal case by assuming that
\begin{equation}                                      \label{beta}
\phi(s)  \le  e^{\a^t - \a}  \quad \text{for $t \ge 1$.}
\end{equation}

\begin{theorem}                               \label{by lorentz}
  Let $\phi$ and $\psi$ be generating functions satisfying
  \eqref{normalization}  and \eqref{beta}.
  Let $\F$ be a class of functions $1$-bounded in $\L_\psi$.
  Then for $0 < t < 1/2$
  $$
  D(F, \L_\phi, t)  \le  C \a \; \vc(F, ct) \cdot \log (\phi | \psi)(t/2)
  $$
\end{theorem}

\remarks 1. No nontrivial estimate is possible when $\phi = \psi$.
Indeed, even in the simplest case when $\Omega$ is finite and
$\mu$ is uniform, let us take $F$ to be the collection of the
functions $f_\omega = \d_\omega / \|\d_\omega\|_{\L_\phi}$,
$\omega \in \Omega$, where $\d_\omega$ is the function that takes
value $1$ at $\omega$ and $0$ elsewhere. Clearly, $F$ is
$1$-bounded in $\L_\phi$ and has combinatorial dimension $d(F,t) =
1$ for any $0 < t < 1$. However, $\|f_\omega -
f_{\omega'}\|_{\L_\phi} \ge 1$ for $\omega \ne \omega'$. Hence
$D(F, \L_\phi, 1/2) = \log |F| = \log |\Omega|$. This can be
arbitrarily large.

2. To see the sharpness of Theorem \ref{by lorentz}, notice that
for {\em some} probability measure $\mu$ on $\Omega$,
$$
D(F, \L_\phi, t)  \ge  c \; \vc(F, Ct).
$$
A simple argument can be found in \cite{T 02} Proposition 1.4.

\qquad

In the extremal case of Theorem \ref{by lorentz}, when $F$ is
$1$-bounded in $L_\infty$, the comparison function becomes just
$\phi(t)$, which gives
\begin{corollary}
  Let $\phi$  be a generating function satisfying \eqref{beta}
  and such that $\phi(1)=1$.
  Let $\F$ be a class of functions $1$-bounded in $L_{\infty}$.
  Then for $0 < t < 1/2$
  $$
  D(F, \L_\phi, t)  \le  C \a \; \vc(F, ct) \cdot \log \phi (t/2).
  $$
\end{corollary}

We use Theorem \ref{by lorentz} for classical Lorentz spaces
$L_{p,\infty} = L_{p,\infty}(\Omega, \mu)$ generated by $\phi(t) = t^p$.

\begin{corollary} \label{weak L_p}
  Let $1 \le p < q \le \infty$.
  Let $\F$ be a class of functions $1$-bounded in $L_{q,\infty}$.
  Then for $0 < t < 1/2$
  $$
  D(F, L_{p,\infty}, t)
  \le  C_{p,q} \; \vc(F,ct) \cdot \log(1/t)
  $$
  where
  \[
  C_{p,q} = C \left( \frac{p^2 q}{q-p} \right).
  \]
\end{corollary}
\proof We apply Theorem \ref{by lorentz} to the functions $\phi
(t)=t^p$ and $\psi(t)=t^q$. In this case the comparison function
becomes $(\phi | \psi)(t) = t^{pq / (p-q)}$. To complete the
proof, notice that \eqref{beta} holds with $\a=p$.
\endproof

Our main interest is in the $L_p$ spaces, for which we obtain
\begin{corollary}               \label{cor Lp-Lq}
  Let $1 \le p < q \le \infty$.
  Let $\F$ be a class of functions $1$-bounded in $L_q$.
  Then for $0 < t < 1/2$
  \begin{equation}                                          \label{Lp-Lq}
    D(F, L_p, t)  \le  C_{p,q} \; \vc(F,c_{p,q}t) \cdot \log(1/c_{p,q}t)
  \end{equation}
  where
  \[
  C_{p,q} = C \left( \frac{p^2 q}{q-p} \right), \ \ \ \
  c_{p,q} = c \min \left( 1, \left(\frac{q-p}{p}\right)^{1/p} \right).
  \]
\end{corollary}
In the next section, this estimate will be applied in an important
partial case, when $p$ is a nontrivial proportion of $q$. In that case,
say if $p \le 0.99 q$, inequality \eqref{Lp-Lq} reads
\begin{equation} \label{L_p estimate}
    D(F, L_p, t)  \le  C p^2 \; \vc(F,ct) \cdot \log(1/t).
\end{equation}
The history of estimates obtained prior to
Corollary~\ref{cor Lp-Lq} and \eqref{L_p estimate} is outlined
in Section~\ref{s:overview} after \eqref{Lp}.

\qquad

\proof Since $\F$ is $1$-bounded in $L_q$, it is also
$1$-bounded in $L_{q,\infty}$.
Let $p'$ be so that $p < p' < q$. Fix an $f \in L_{p', \infty}$ with
$\|f\|_{p',\infty} \le 1$. Then
\begin{align*}
  \|f\|_p^p &\le \int_{\{\w :\; |f(\w)| \le 1\}} |f(\w)|^p d \mu +
  \int_1^{\infty} p t^{p-1} \mu \{\w :\; |f(\w)| \ge t\} dt \\
  &\le 1+ p \int_1^{\infty}t^{p-1-p'} dt \le \frac{p'}{p'-p}.
\end{align*}
Taking the $p$-th root we conclude that if $f,g$ are $t$-separated
in $L_p$, then they are $(b_{p,p'} t)$-separated in $L_{p',\infty}$,
where
$$
b_{p,p'} = \left( \frac{p'-p}{p'} \right)^{1/p}.
$$
Thus $D(F, L_p, t)  \le D(F, L_{p',\infty}, b_{p,p'} t)$.
Then the application of Corollary \ref{weak L_p} with $p'$ and $q$ gives
$$
D(F, L_p, t)  \le  B_{p',q} \; \vc(F,b_{p,p'}t) \cdot \log(1/b_{p,p'}t)
$$
with
$$
B_{p',q} = C \left( \frac{{p'}^2 q}{q-p} \right).
$$
If we choose $p' = \min(2p, \frac{p+q}{2})$ then a direct check shows that
$$
B_{p',q} \le C \left( \frac{12 p^2 q}{q-p} \right), \ \ \
b_{p,p'} \ge  \min \left( \frac{1}{2}, \left( \frac{q-p}{4p} \right)^{1/p}
                    \right).
$$
This completes the proof.
\endproof

To prove Theorem \ref{by lorentz}, we first reduce the size of the
domain $\Omega$ (which can be assumed finite) by means of a
probabilistic selection and then apply the covering Theorem
\ref{by tower}.

In the probabilistic selection, we use a standard independent
model. Given a finite set $I$ and a parameter $0 < \d < 1$, we
consider selectors $\d_i$, $i \in I$, which are independent
$\{0,1\}$-valued random variables with $\E \d_i = \d$. Then the
set $J=\{i \in I : \d_i = 1\}$ is a random subset of $I$ and its
average cardinality is $s=\d |I|$. We call $J$ a random set of
expected cardinality $s$.

\begin{lemma}                              \label{binomial}
  Let $0 < \e \le 1$. For $t \ge \e \d m$,
  $$
  \P \Big\{ \Big| \sum_{i=1}^m (\d_i - \d) \Big| > t \Big\}
  \le  2 \exp(-c \e t)
  $$
  where $c > 0$ is an absolute constant.
\end{lemma}

\proof This follows from Prokhorov-Bennett inequality.
Let $(X_i)$ be a finite sequence of real valued independent
mean zero random variables such that $\|X_i\|_\infty \le a$
for every $i$. If $b^2 = \sum_i \E X_i^2$, then for all $t > 0$
\begin{equation}            \label{prokhorov-bennett}
  p := \P \Big\{ \sum_i X_i > t \Big\}
  \le \exp \Big[ t/a - (t/a + b^2/a^2) \log(1 + at/b^2) \Big]
\end{equation}
which is less than $\exp(-t^2 / 4b^2)$ if $t \le b^2/2a$
(see e.g. \cite{LT} 6.3).

We apply Prokhorov-Bennett inequality for $X_i = \d_i - \d$
and with $a = 1$, $b^2 = \d m$.
Consider two cases:

1) $\e \d m \le t \le 8 \d m$. Since in that case $t/16 \le b^2/2a$,
we have
$$
p  \le  \P \Big\{ \sum_i X_i > t/16 \Big\}
\le  \exp(-t^2 / 64 \d m)
\le  \exp(-\e t / 64)
$$
because $t \ge \e \d m$.

2) $t > 8 \d m$. Then $\log(1 + at/b^2) = \log(1 + t / \d m) > 2$, hence
$$
p  \le  \exp \Big[ -(t/a) \Big( \log(1 + at/b^2) - 1 \Big) \Big]
<  \exp(-t).
$$

Thus for all $t > \e \d m$ we have $p \le \exp(-c \e t)$.
Repeating the argument for $-X_i$, we conclude the proof.
\endproof

\begin{lemma}                                   \label{discrepancy}
   There exist absolute constants $C, c > 0$ for which
   the following holds.
   Let $\g > 0$ and let $\QQ$ be a system of subsets
   of $\{1, \ldots, n\}$ such that
   $$
   |S|  \ge  \g n \ \ \
   \text{for all $S \in \QQ$.}
   $$
   If $\s$ is a random subset of $\{1,\ldots,n\}$
   of expected cardinality $k$ satisfying
   $|\QQ| \le 0.001 \cdot \exp(c \g k)$,
   then with probability at least $0.99$ we have
   $$
   \frac{|S \cap \s|}{|\s|}
   \ge  0.99 \, \frac{|S|}{n}
   \ \ \ \text{for all $S \in \QQ$}.
   $$
\end{lemma}

\proof Let $0 < \d < 1/2$ and set $\d_1, \ldots, \d_n$ to be
$\{0,1\}$-valued independent random variables with $\E \d_i = \d$
for all $i$. Let $\d = k/n$; consider the random set $\s = \{ i:
\; \d_i = 1 \}$. For any set $S  \subset \{1, \ldots, n \}$, $|S
\cap \s|  = \sum_{i \in S} \d_i$. By Lemma \ref{binomial} applied
to a sum over $S$ instead of $\{1,\ldots, m\}$ and with $t = 0.001
\d |\QQ|$, there is an absolute constant $c_0 > 0$ such that
$$
\P \{  |S \cap \s| < 0.999 \d |S| \}
   \le   2 \exp (- c_0 \d |S|).
$$
Since for every $S \in \QQ$, $|S| \ge \g n$, then
$$
\P \Big\{  |S \cap \s|  <  0.999 \frac{k}{n} \cdot |S| \Big\}
   \le   2 \exp (- c_0 \g k).
$$
Therefore
$$
\P \Big\{  \forall S \in \QQ, \; |S \cap \s| \ge  0.999
\frac{k}{n} \cdot |S| \Big\}
   \ge  1 - 2 |\QQ| \exp (- c_0 \g k)
    \ge 0.998
$$
provided $c \le c_0 / 2$. Also, $|\s| \le 1.001 k$ with
probability at least $0.999$ since $k$ can be assumed sufficiently
large. This completes the proof.
\endproof

Given a fininte set $I$, we will now work with Lorentz spaces
$\L_\phi(I) = \L_\phi(I,\mu)$, where $\mu$ is the uniform measure
on $I$. The following two lemmas reduce the size of $I$
while keeping both the boundedness of the class $F$ in the $\L_\psi$-norm
and the separation of $F$ in the $\L_\phi$-norm.

\begin{lemma}                                 \label{one function}
  Let $\psi$ be a generating function.
  Let $f$ be a function on a finite set $I$ such that
  $$
  \|f\|_{\L_\psi(I)} \le 1.
  $$
  If $\s$ is a random subset of $I$ of expected
  cardinality $k > C$, then
  with probability at least $0.9$ we have
  $$
  \|f\|_{\L_\psi(\s)} \le C.
  $$
\end{lemma}

\proof
Let $a > 2$ be a parameter to be chosen later and
let $\d, \d_j$ be as in the proof of Lemma \ref{discrepancy}.
For $s \in \Z$ define the set
\[
I_s=\{ j \in I \; : \; |f(j)| > 2^s \}.
\]
Since $\|f\|_{\L_\psi(I)} \le 1$, we have
\begin{equation}                                \label{Is small}
|I_s|  \le  \frac{n}{\psi(2^s)}.
\end{equation}
Define also the event $A_s$ as
$$
A_s = \Big\{ |I_s \cap \s| > \frac{a \d n}{\psi(2^s)} \Big\}.
$$
We want to bound the probability that at least one $A_s$ occurs.
Let $r$ be the maximal number such that $\d |I_r| \ge 0.01$. Then
\[
\P \{ \forall j \in I_{r+1} \; \d_j=0 \} =(1-\d)^{|I_{r+1}|}
\ge e^{-\d |I_{r+1}|} \ge e^{-0.01} > 0.99.
\]
If $A_s$ occurs for some $s > r$ then $I_s \cap \s$ is nonempty,
hence the larger set $I_{r+1} \cap \s$ is nonempty, which happens
with probability at most $0.01$.
Thus
\begin{equation}                                \label{s < r}
\P \Big( \bigcup_{s > r} A_s \Big) \le 0.01.
\end{equation}
To bound $\P(A_s)$ with $s \le r$, we will apply Lemma \ref{binomial}
with $m = |I_s|$ and $t = \frac{a \d n}{2 \psi(2^s)}$. Note that
$$
\P (A_s) = \P \Big\{ \sum_{j \in I_s} (\d_j-\d)
                     > \frac{a \d n}{\psi(2^s)} - \d m \Big\}.
$$
By \eqref{Is small}, we have
$\frac{a \d n}{\psi(2^s)} - \d m  \ge  t  \ge  \d m$.
Then Lemma \ref{binomial} gives
$$
\P (A_s) \le \exp \Big( -\frac{c a \d n}{\psi(2^s)} \Big).
$$
By the convexity of $\psi$,
\begin{equation}                                   \label{psi regular}
\psi(wx) \ge w \psi(x) \ \ \ \text{for all $x \ge 0$ and $w \ge
1$.}
\end{equation}
Thus $\psi(2^s) \le 2^{s-r} \psi(2^r)$.
Then using \eqref{Is small} and the fact that $\d |I_r| \ge 0.01$,
we obtain
$$
\P (A_s) \le \exp \Big( - 2^{r-s} \frac{c a \d n}{\psi(2^r)} \Big)
  \le  \exp ( - 2^{r-s} c a \d |I_r| )
  \le  \exp ( - 0.01 c a 2^{r-s} ).
$$
So if $a$ is taken large enough then
$\sum_{s=-\infty}^r \P(A_s) \le 0.01$. Combining this
with \eqref{s < r}, we conclude that
$$
\P \Big( \bigcup_{s \in \Z} A_s \Big) \le 0.02.
$$
In addition, by Lemma \ref{binomial} we have
$\P \{ |\s| < \frac{1}{2} \d n \}  \le 0.02$ since $k$ can
be assumed large enough.

Now suppose that $|\s| \ge \frac{1}{2} \d n$
and that none of the events $A_s$ occur,
which happens with probability at least
$1-0.02-0.02=0.96$. Fix any $t > 0$ and find an integer $s$
so that $2^s \le t < 2^{s+1}$. By the definitions
of $A_s$, $I_s$ and by \eqref{psi regular},
$$
| \{i \in \s : |f(i)| > t \} |
\le |I_s \cap \s|
\le \frac{a \d n}{\psi(2^s)}
\le \frac{2 a |\s|}{\psi(2^s)}
\le \frac{2 a |\s|}{\psi(t/2)}
\le \frac{|\s|}{\psi(t/4a)}.
$$
This means that
$\| f\|_{\L_{\psi}(\s)} \le 4a$.
\endproof

\qquad

\begin{lemma}                                   \label{reduction}
  Let $\phi$, $\psi$ be Lorentz functions.
  Let $\F$ be a class of functions on a finite set $I$,
  which is $1$-bounded in the $\L_\psi(I)$ norm.
  Assume that
  \begin{equation}                                \label{x large}
    \|x\|_{\L_\phi(I)} \ge t   \ \ \
    \text{for all $x \in F$.}
  \end{equation}
  If $\s$ is a random subset of $I$ of expected cardinality $k$
  satisfying $|F| \le 0.001\exp \big( \frac{ck}{(\phi|\psi)(t)} \big)$,
  then with probability at least $0.99$ we have
  $$
    \|x\|_{\L_\phi(\s)}  \ge  0.99 \|x\|_{\L_\phi(I)}  \ \ \
    \text{for all $x \in F$.}
  $$
\end{lemma}
This lemma  will be applied to the difference set $A - A$ of a
$t$-net $A$ of the class $F$ in the theorem.

\qquad

\proof We can assume that $I = \{1, \ldots, n\}$. Fix an $x \in
\F$. Since $|x| / \|x\|_{\L_\phi(I)} = 1$, there exists an $s =
s(x) > 0$ such that
\begin{equation}                                            \label{mux large}
\mu \Big\{ \frac{|x|}{\|x\|_{\L_\phi(I)}}  > s \Big\}
  \ge  \frac{1}{\phi(s)}.
\end{equation}
On the other hand, since $\|x\|_{\L_\phi(I)} \ge t$ and
$\|x\|_{\L_\psi(I)} \le 1$, the measure in \eqref{mux large}
is majorized by
$$
\mu \{ |x| > t s \}  \le  \frac{1}{\psi(t s)}.
$$
Hence $\phi(s) \ge \psi(ts)$ and therefore
\begin{equation}                                            \label{phi small}
  \phi(s) \le (\phi|\psi)(t).
\end{equation}
Now consider the family of subsets of $I$ defined as
$$
S(x)  =  \Big\{ i : \frac{|x(i)|}{\|x\|_{\L_\phi(I)}}  > s(x)
\Big\}, \ \ \ x \in \F.
$$
By \eqref{mux large} and \eqref{phi small},
for every $x \in F$
$$
|S(x)|  \ge  \frac{n}{\phi(s(x))} \ge \frac{n}{(\phi|\psi)(t)}.
$$
Let $\mu_\s$ denote the uniform probability measure on $\s$. Lemma
\ref{discrepancy} implies that whenever $|F| \le 0.001\exp \big(
\frac{ck}{(\phi|\psi)(t)} \big)$, a random subset $\s$ of $I$ of
average cardinality $k$ satisfies with probability at least $0.99$
that
\begin{align*}
\mu_\s \Big\{ \frac{|x|}{\|x\|_{\L_\phi(I)}}  > s(x) \Big\}
  &=    \frac{|S(x) \cap \s|}{|\s|}
   \ge  0.99 \frac{|S(x)|}{n}\\
  &\ge  0.99 \frac{1}{\phi(s(x))}
   \ge  \frac{1}{\phi(s(x)/0.99)} \ \ \
\text{for all $x \in \F$.}
\end{align*}
Hence
$$
\|x\|_{\L_\phi(\s)}  \ge  0.99 \|x\|_{\L_\phi(I)}  \ \ \
\text{for all $x \in F$.}
$$
The proof is complete.
\endproof

\qquad

\noindent {\bf Proof of Theorem \ref{by lorentz}.} We may assume
that $\Omega$ is finite. By splitting the atoms
of $\Omega$ (by replacing an atom $\w$ by, say, two atoms $\w_1$
and $\w_2$, each carrying measure $\frac{1}{2} \mu(\w)$ and by
defining $f(\w_1) = f(\w_2) = f(\w)$ for $f \in F$), we can make
the measure $\mu$ almost uniform without changing neither the
covering numbers nor the combinatorial dimension of $F$. So, we
can assume that $\mu$ is the uniform measure on $\Omega$.

Let $A$ be a $t$-separated subset of $F$ (which means that
$\|f-g\|_{\L_\phi(\Omega)} \ge t$ for all $f \ne g$ in $F$) of
size
$$
\log|A| = D(F, \L_\phi(\Omega), t).
$$
 The
difference set $\frac{1}{2} (A - A) \setminus \{0\} = \{
\frac{1}{2} (f-g) : f \ne g; \, f,g \in A \}$ satisfies the
assumptions of Lemma \ref{reduction} with $t/2$ in \eqref{x
large}. Then for $k$ defined by
\begin{equation}                                \label{size of A}
  |A|^2  =  0.001 \exp \Big( \frac{ck}{(\phi|\psi)(t/2)} \Big),
\end{equation}
a random subset $\s$ of $\Omega$ of average cardinality $k$
satisfies with probability at least $0.99$ that
$$
  \big\| \frac{1}{2}(f-g) \big\|_{\L_\phi(\s)}
  \ge  0.99 \big\| \frac{1}{2}(f-g) \big\|_{\L_\phi(\Omega)}
  \ge  0.99 \frac{t}{2}
  \ge  \frac{t}{3} \ \ \ \text{for all $f \ne g$ in $A$.}
$$
This means that
\begin{equation}                                 \label{A separated}
  \text{$A$ is $\frac{t}{3}$-separated in $\L_\phi(\s)$}
\end{equation}
and in particular
$$
D(F, \L_\phi(\s), t/3 ) \ge \log|A| = D(F, \L_\phi(\Omega), t).
$$
The advantage of the left hand side is that the size of $\s$ is
controlled via \eqref{size of A}.

We need also to keep $A$ well bounded in $\L_\psi(\sigma)$. Denote
by $\E_\s$ the average over the random set $\s$, that is over the
selectors $\d_i$. By Lemma \ref{one function}
\[
  \E\; |\{ f \in A \;: \|f\|_{\L_\psi(\s)} \le C\}|
  =\sum_{f \in A} \P (\|f\|_{\L_\psi(\s)} \le C) \ge 0.9|A|.
\]
Therefore with probability at least $0.8$,
\begin{equation}                \label{half small}
\text{at least a half of the functions in $A$ have norm
      $\|f\|_{L_\psi(\sigma)} \le C$.}
\end{equation}
Since $k/2 \le |\s| \le 2k$ holds with probability at least $0.9$,
there exists a realization of $\s$ that satisfies simultaneously
this property, \eqref{A separated} and \eqref{half small}.
Let $B$ be the set consisting of $\frac{6}{t} f$,
where $f$ are the functions satisfying \eqref{half small}.

Summarizing, there exists a subset $\s$ of $\Omega$ and a set
$B$ such that

\begin{itemize}
  \item{$B$ is a subset of $\frac{6}{t} A$,}
  \item{$B$ is $(C/t)$-bounded in $L_\psi(\sigma)$,}
  \item{$B$ is $2$-separated in $\L_\phi(\sigma)$,}
  \item{$|B| \ge |A|/2 \ge c' \exp \big( \frac{c|\s|}{2(\phi|\psi)(t/2)}
        \big)$.}
\end{itemize}

We can clearly assume that $\s = \{1, \ldots, n\}$ and
realize the space $\L_\phi^n = \L_\phi(\{1, \ldots, n\})$ as
$\R^n$ equipped with the Lorentz norm $\L_\phi$.
Applying the covering Theorem \ref{by tower}, we get
\begin{equation}                                \label{B covered}
  N(B, \tower^{\a})  \le  \Sigma(B)^{\a}.
\end{equation}
Since $B$ is $2$-separated in $\L_\phi^n$ and by \eqref{beta} the
norm in this space is bounded by the $\tower^\a$ norm, the set $B$ is
also $2$-separated in the $\tower^\a$ norm. Hence
\begin{equation}                \label{N=B}
N(B, \tower^{\a}) = |B|.
\end{equation}
The right hand side of \eqref{B covered} can be estimated through
Lemma \ref{A in L1}.
By \eqref{normalization} and convexity, $\psi(t) \ge t$ for $t \ge 1$.
Thus $C \|f\|_{L_\psi} \ge \|f\|_{L_{1/2}}$ for all functions $f$.
Hence
$$
B \subseteq (C/t) \, \Ball(L_\psi^n)
  \subseteq (C'/t) \, \Ball(L_{1/2}^n).
$$
Hence by \eqref{B covered}, \eqref{N=B} and Lemma \ref{A in L1},
\begin{equation}                \label{B controlled}
  |B| \le  \Big( \frac{C n}{t v} \Big)^{3 v}, \ \ \
  \text{where $v = \vc(B) = \vc(B,1) \le \vc(A, t/6)$.}
\end{equation}
We also have a lower bound $|B| \ge c' \exp(a n)$ with $a =
\frac{c/2}{(\phi|\psi)(t/2)}$. Taking logarithms of the upper
and the lower bounds, we obtain $a n/v  \le  C \log(C n/tv)$,
from which it follows that
$$
\frac{n}{v}  \le  \frac{C}{a} \log \Big(\frac{C}{t a}\Big).
$$
Plugging this back into \eqref{B controlled}, we obtain
$$
|B|  \le  \Big[ \Big(\frac{C}{t a}\Big)
           \log \Big(\frac{C}{t a}\Big) \Big]^{C v}
\le  \Big(\frac{C}{t a}\Big)^{C_1 v}.
$$
Note that
$$
(\phi|\psi)(t)  \ge  1/t \ \ \ \text{for $0 < t < 1$}.
$$
Indeed, $\phi(\frac{1}{t}) \ge 1 = \psi(t \cdot \frac{1}{t})$
which implies $(\phi|\psi)(t) \ge \phi(\frac{1}{t}) \ge
\frac{1}{t}$.

Therefore $t \ge (2/c) a$ and finally
$$
|B| \le a^{-C v} \le (\phi|\psi)(t/2)^{C v}.
$$
On the other hand, by the construction
$$
\log|B| \ge \log \big(\frac{1}{2} |A|\big)
\ge c D(F, \L_\phi(\Omega), t).
$$
This completes the proof.
\endproof

\section{Random processes and the uniform entropy}
       \label{s:gaussian}

Here we prove our main results, Theorems \ref{i:integral=integral} and
\ref{i:entropy=dimension}, which compare the uniform entropy $D(F,t)$
to the combinatorial dimension $\vc(F,t)$. One direction of this comparison
is easy: for every class of functions $F$ and every $t>0$,
\begin{equation}            \label{entropy>dimension}
D(F,t)  \ge  c \; \vc(F, 2t)
\end{equation}
where $c>0$ is an absolute constant, see \cite{T 02}.

The reverse inequality is not true in general even for $\{0,1\}$
classes. Let, for example, $F$ be the collection of $n$
characteristic functions
$\one_{ \{ i \} }$ of the singletons $i \in \{1,\ldots, n\}$.
Then for $0 < t < n^{-1/2}$ we have $D(F,t) = \log n$ while
$\vc(F,t) = 1$.

Nevertheless, we are able to show that the reverse to
\eqref{entropy>dimension} holds:
1) under a minimal regularity of $F$, and
2) always after taking integrals on both sides.

\paragraph{Integral equivalence}
The following is a general form of Theorem \ref{i:integral=integral}.

\begin{theorem}             \label{integral=integral}
  For every class $F$ and for any $b \ge 0$,
  \begin{equation}          \label{eq integral=integral}
    \int_b^\infty \sqrt{D(F,t)} \; dt
    \le C \int_{cb}^\infty \sqrt{\vc(F,t)} \; dt.
  \end{equation}
\end{theorem}

The proof of Theorem \ref{integral=integral} is based on the
following
\begin{lemma} \label{entropy iteration}
Let $a>2$ and let $F$ be a function class. Then for all $t>0$:
\begin{equation}                \label{D via sum}
D(F,t)  \le  C \log a \sum_{j=0}^\infty 4^j \vc(F, c a^j t).
\end{equation}
\end{lemma}
\proof The proof of the Lemma uses an iteration argument. It
relies on the fact valid for arbitrary sets $K$, $D$ and $L$ in
$\R^n$:
\begin{equation}                    \label{double cover}
  N(K,D)  \le  N(K,L) \sup_{z \in \R^n} N( (K+z) \cap L, D).
\end{equation}
To check this, first cover $K$ by translates of $L$ and then cover
the intersection of $K$ with each translate by appropriate translates of $D$.

It will be easier to work with the ``covering'' analog of
$D(F,t)$, so we define a covering version of $D(F,X,t)$ in
\eqref{D(F,X,t)} as
$$
D'(F,X,t)
= \log \sup \Big( n \mid \exists f_1, \ldots, f_n \in X \
\forall f \in F \ \exists i \  \|f-f_i\|_X \le t \Big).
$$
By \eqref{covering vs packing},
\begin{equation}                    \label{D vs D'}
D(F,X,2t)  \le  D'(F,X,t)  \le  D(F,X,t).
\end{equation}

We can clearly assume the domain $\Omega$ to be finite. Fix the
underlying probability $\mu$ on $\Omega$ and $t >0$. For $j = 1,2,
\ldots$ define
$$
t_j = a^{j-1} t  \ \ \ \text{and} \ \ \ X_j = L_{2^j} (\Omega,\mu).
$$
We estimate $D(F,L_2(\Omega,\mu),t)$ by \eqref{double cover}:
$$
D'(F,X_1,t_1)
\le  D'(F,X_2,t_2) + \sup_h D' \big( (F+h) \cap t_2 \Ball(X_2), X_1, t_1\big)
$$
where the supremum is over all functions $h$ on the (finite) domain $\Omega$.
Iterating this inequality, we obtain
\begin{equation}                    \label{iteration}
D'(F,X_1,t_1)  \le  \sum_{j=1}^\infty \sup_h D'(F_j(h), X_j, t_j)
\end{equation}
where $F_j(h) = (F+h) \cap t_{j+1} \Ball(X_{j+1})$.
Obviously, the class $F_j(h)$ is $t_{j+1}$-bounded in $X_{j+1}$.
Then applying \eqref{L_p estimate} to the class $t_{j+1}^{-1} F_j(h)$
with $p = 2^j$ and $q = 2^{j+1}$, we obtain
$$
D'(F_j(h), X_j, t_j)
\le  C 4^j \vc(F_j(h), c t_j) \cdot \log(t_{j+1}/t_j)
\le  C 4^j \vc(F, c t_j) \cdot \log a
$$
because apparently $\vc(F_j(h), s) \le \vc(F, s)$ for every $s$.
To complete the proof we substitute the previous inequality into
\eqref{iteration} and use \eqref{D vs D'}.
\endproof

\qquad

\noindent {\bf Proof of Theorem \ref{integral=integral}.\ }
Applying Lemma \ref{entropy iteration} and Jensen's inequality, we
have for $a=3$:
\begin{align*}
\int_b^\infty \sqrt{D(F,t)} \; dt
  &\le  C \sqrt{\log a} \;\sum_{j=0}^\infty
        2^j \int_b^\infty \sqrt{\vc(F, c a^j t)} \;dt  \\
  &=    C \sqrt{\log a} \;\sum_{j=0}^\infty
        (2/a)^j \int_{c a^j b}^\infty \sqrt{\vc(F, u)} \;du
  \le C \int_{cb}^\infty \sqrt{\vc(F,u)} \; du.
\end{align*}
\endproof

\remark
The square roots in \eqref{eq integral=integral} can of course
be replaced by any other equal positive powers. The present
form of \eqref{eq integral=integral} was chosen to match
Dudley's entropy integral, see Theorem \ref{supremum} below.

\qquad

Theorem \ref{i:donsker} follows from Theorem
\ref{integral=integral} as explained in the introduction. The gap
between the sufficient and the necessary conditions in Theorem
\ref{i:donsker} is known to be needed in general (at least in
the description of universal Donsker classes, see \cite{Du 99}
Propositions 10.1.8 and 10.1.14).

\paragraph{Pointwise equivalence}
A similar argument, which we give now, completes the proof
of the other main result of the paper, Theorem \ref{i:entropy=dimension}.
Although \eqref{entropy>dimension} can not be reversed in general,
a remarkable fact is that it can be reversed if the combinatorial
dimension is polynomial in $t$.

\begin{theorem}             \label{polynomial}
  Let $F$ be a class of functions and $t > 0$.
  Assume that there exist positive numbers $v$ and $\a \le 1$
  such that
  \begin{equation}          \label{polynomial assumption}
    \vc(F, tx)  \le  v x^{-\a} \ \ \
    \text{for all $x \ge 1$.}
  \end{equation}
  Then
  $$
  D(F, Ct)  \le  (C/\a) \, v.
  $$
\end{theorem}

\proof Applying Lemma \ref{entropy iteration} with $a = 5^{1/\a}$
and estimating the combinatorial dimension via \eqref{polynomial
assumption}, we have
$$
\vc(F, t/c)  \le  C \log a \sum_{j=0}^\infty 4^j v a^{-\a j}
\le (C/\a) \, v.
$$
This completes the proof.
\endproof

The following is a general form of Theorem \ref{i:entropy=dimension}.
It improves upon Talagrand's inequality proved in \cite{T 87},
see \cite{T 02}.

\begin{corollary}
  Let $F$ be a class of functions and $t > 0$.
  Assume that there exists a decreasing function $v(t)$ and a number $a > 2$
  such that
  \begin{equation}          \label{regularity}
  \vc(F,s) \le v(s)  \ \ \ \text{and} \ \ \
  \vc(as)  \le  \frac{1}{2} \, \vc(s)
  \ \ \ \text{for all $s \ge t$.}
  \end{equation}
  Then
  $$
  D(F, Ct)  \le  C \log a \cdot v(t).
  $$
\end{corollary}

\proof Applying \eqref{regularity} recursively, we have $v(a^j t)
\le  \frac{1}{2^j} \, v(t)$ for all $j = 0,1,2,\ldots$. Let $x \ge
1$ and choose $j$ so that $a^j \le x \le a^{j+1}$. Then
$$
  2^{-j} \ge x^{-\frac{\log 2}{\log a}} \ge 2^{-j-1},
$$
so,
$$
\vc(F, tx)  \le  v(tx) \le v(ta^j)  \le C v(t) \cdot 2^{-j} \le 2C
\, v(t) \cdot x^{-\frac{\log 2}{\log a}}.
$$
 The conclusion follows by Theorem
\ref{polynomial}.
\endproof

\paragraph{A combinatorial bound on Gaussian processes}
A quantitative version of Theorem \ref{i:donsker} is the following
bound on Gaussian processes indexed by $F$ in terms of the
combinatorial dimension of $F$.

Let $F$ be a class of functions on an $n$-point set $I$.
The standard Gaussian process indexed by $f \in F$ is
$$
X_f  =  \sum_{i \in I} g_i f(i)
$$
where $(g_i)$ are independent $N(0,1)$ random variables. The
problem is to bound the supremum of the process $(X_f)$ normalized
by the standard deviation as in the Central Limit Theorem:
$$
E(F) = n^{-1/2} \; \E \sup_{f \in F} X_f.
$$

\begin{theorem}             \label{supremum}
  For every class $F$,
  \begin{equation}          \label{eq supremum}
  E(F)  \le  C \int_0^\infty \sqrt{\vc(F,t)} \; dt
  \end{equation}
  where $C$ is an absolute constant.
  Moreover, $0$ can be replaced by $c n^{-1/2} E(F)$,
  where $c>0$ is an absolute constant.
\end{theorem}

\proof
By Dudley's entropy integral inequality,
\begin{equation}                \label{dudley}
E(F)  \le  C \int_{n^{-1/2} E(F)}^\infty
           \sqrt{D(F,L_2(\mu),t)} \; dt
\end{equation}
where $\mu$ is the uniform probability measure on $I$,
see \cite{MV}.
Then the proof is completed by Theorem \ref{integral=integral}.
\endproof

In 1992, M.~Talagrand proved Theorem \ref{supremum}
for uniformly bounded convex classes
and up to an additional factor of $\log^M(1/t)$ in the integrand;
this was a main result of \cite{T 92}.
The absolute constant $M$ was reduced to $1/2$ in \cite{MV}.
Theorem \ref{supremum} is optimal.
We emphasize its important meaning:
\begin{quote}
  {\em In the classical Dudley's entropy integral, the entropy can
  be replaced by the combinatorial dimension.}
\end{quote}

\paragraph{Optimality of the bound on Gaussian processes}
We conclude this section by showing the sharpness of
Theorem \ref{supremum}. For every $n$ one easily finds
a class $F$ for which the inequality in \eqref{supremum}
can be reversed -- this is true e.g. for $F=\{-1,1\}^I$.
More importantly, the integral in \eqref{supremum} can not be
improved in general to the (Sudakov-type) supremum
$\sup t \sqrt{\vc(F,t)}$. This is so even if we replace the
Gaussian process $X_f$ by the Rademacher process
$$
Y_f  =  \sum_{i \in I} \e_i f(i)
$$
where $(\e_i)$ are independent symmetric $\pm 1$ valued random
variables. The average supremum of such process,
$$
E_\rad (F) = n^{-1/2} \; \E \sup_{f \in F} Y_f,
$$
is well known to be majorized by that of the Gaussian process:
$E_\rad (F) \le CE(F)$ (see \cite{LT} Lemma 4.5).

\begin{proposition}         \label{no sudakov}
  For every $n$, there exists a class $F$ of functions
  on $\{1,\ldots,n\}$ uniformly bounded by $1$ and such that
  \begin{equation}          \label{eq no sudakov}
    E_\rad  (F)  \ge  c_1 \int_0^\infty \sqrt{\vc(F,t)} \; dt
    \ge c \log n \cdot \sup_{t > 0} t \sqrt{\vc(F,t)}.
  \end{equation}
\end{proposition}
Our example will be constructed as sums of random vertices of
the discrete cube with quickly decreasing weights.

We shall bound $E_\rad (F)$ from below via a Sudakov type
minoration for Rademacher processes.
Let $D = \Ball(L_2^n) = \sqrt{n} B_2^n$.
Proposition 4.13 of \cite{LT} with $\e = \frac{1}{2} n$
states the following:

\begin{fact}  \label{Rademacher Sudakov}
  If $A$ is a subset of $\R^n$ and
  \begin{equation}  \label{A in the cube}
    \sup_{x \in A} \|x\|_{\infty} \le c_1 \frac{\sqrt{n}}{E_\rad(A)},
  \end{equation}
  then
  \begin{equation}              \label{rad entropy}
    \sqrt{\log N \big( A, \frac{1}{2} D \big)} \le CE_\rad(A).
  \end{equation}
\end{fact}

The entropy in \eqref{rad entropy} will be estimated
in a standard way:

\begin{fact}  \label{entropy of random points}
  There exists an absolute constant $\a$ such that 
  the following holds. 
  Let $A$ be a set of $N \le e^{\a n}$ 
  random vertices of the discrete cube
  $\{-1,1\}^n$, i.e. $A$ consists of $N$ independent copies
  of a random vector $(\e_1, \ldots, \e_n)$. Then with
  probability at least $1/2$,
  $$
  N(A, \frac{1}{2} D)  \ge  \sqrt{N}.
  $$
\end{fact}

\proof Obviously, we can assume that $N \ge 2$. Assume that the
event
\begin{equation}          \label{entropy small}
    N(A, \frac{1}{2} D) \le \sqrt{N}
\end{equation}
occurs. Then there exists a translate
$D_x' = \frac{1}{2} D + x$ of $\frac{1}{2} D$ which
contains at least $N / N(A, \frac{1}{2} D) \ge \sqrt{N}$ points
from $A$. Set $A'=A \cap D_x'$. By dividing the set $A'$ into
pairs in an arbitrary way, we can find a set $\PP$ of $M \ge
\frac{\sqrt{N} - 1}{2}$ pairs $(x,y) \in A' \times A'$, $x \ne y$,
so that each point from $A$ belongs to at most one pair in $\PP$.
Since $A'$ lies in a single translate of $\frac{1}{2} D$, we have
\begin{equation}                \label{distances small}
  \|x-y\|_{L_2^n}  \le 1 \ \ \ \text{for all $(x,y) \in \PP$.}
\end{equation}
Thus if \eqref{entropy small} occurs, then \eqref{distances small}
occurs for some $M$-element set $\PP \subset A \times A$.

Let now $\PP$ be a {\em fixed} set of $M$ disjoint pairs of
elements of $A$. Then
\begin{equation}            \label{prob distances small}
  \P(\text{event \eqref{distances small} occurs})
  = ( \P (\|x-y\|_{L_2^n} \le 1) )^M
\end{equation}
where $x$ and $y$ are independent random vertices of the
discrete cube. Here we used the fact that the pairs in $\PP$
are disjoint from each other and, consequently, are jointly independent.
The probability in \eqref{prob distances small} is easily estimated
using Prokhorov-Bennett inequality \eqref{prokhorov-bennett}:
$$
\P (\|x-y\|_{L_2^n} \le 1)
=  \P \Big( \sum_{i=1}^n |\e_{1i} - \e_{2i}|^2 \le n \Big)
\le e^{-c_1 n}
$$
where $(\e_{1i})$ and $(\e_{2i})$ are independent copies of
the random vector $(\e_i)$.

To estimate the probability that the event \eqref{entropy small}
occurs, note that there is less than $\binom{N^2}{M}$ ways to
choose $\PP$. Therefore
$$
p:= \P \Big( N(A, \frac{1}{2} D) \le \sqrt{N} \Big)
\le \binom{N^2}{M} (e^{-c_1 n})^M
\le \Big( \frac{e N^2}{M} e^{-c_1 n} \Big)^M
\le \big( N e^{-c_2 n} \big)^{3M/2}.
$$
Since $N \ge 2$, we have $M \ge 1$. If $\a \le c_2/2$, we conclude
that $p \le \exp( - \frac{c_2 n}{2} \cdot 3M/2) < 1/2$. This
completes the proof.
\endproof

\begin{corollary} \label{random vertices}
  There exists an absolute constant $\a$ such that 
  the following holds. 
  Let $A$ be a set of $N \le e^{\a n}$ random vertices 
  of the discrete cube.
  Then with probability at least $1/2$,
  $$
  c \sqrt{\log N} \le E_\rad (A) \le C \sqrt{\log N}.
  $$
\end{corollary}

\proof
Since $A \subset \{-1,1\}^n \subset \sqrt{n} B_2^n$,
we have
\begin{equation}                \label{Erad small}
  E_\rad (A) \le C E(A) \le C_1  \sqrt{\log N},
\end{equation}
see \cite{LT} (3.13). To prove the reverse inequality, assume that
$N \le e^{\b n}$ for $\b = \min (\a, (c_1/C_1)^2)$, where $c',
c_1$ and $C_1$ are the absolute constants from Fact  \ref{entropy
of random points}, \eqref{A in the cube} and \eqref{Erad small}
respectively. Then \eqref{Erad small} implies that \eqref{A in the
cube} is satisfied. Hence by Facts \ref{Rademacher Sudakov} and
\ref{entropy of random points}, with probability at least $1/2$ we
have
$$
E_\rad (A) \ge (1/C) \sqrt{\log N(A, \frac{1}{2} D)}
\ge (c / \sqrt{2}) \sqrt{\log N}.
$$
This completes the proof.
\endproof

\qquad

\noindent {\bf Proof of Proposition \ref{no sudakov}. } Fix a
positive integer $n$. Let $k_1$ be the maximal integer so that
$2^{4^{k_1}} \le e^{\a n}$, where $\a$ is an absolute constant
from Corollary \ref{random vertices}. For each $1 \le k \le k_1$
define a set $F_k$ in $\R^n$ as follows. Let $N(k) = 2^{4^k}$ and
let $A_k= \{x_1^k, \ldots, x_{N(k)}^k\}$ be a family of points in
the discrete cube $\{-1,1\}^n$ in $\R^n$ satisfying Corollary
\ref{random vertices}.
 Set
\[
F_k=2^{-k} \cdot A_k \ \ \ \text{and put} \ \ \ F =
\sum_{k=1}^{k_1} F_k
\]
where the sum is the Minkowski sum: $A+B = \{a+b \;:\; a \in A, b
\in B\}$. Then $F$ is a uniformly bounded class of functions on
$\{1,\ldots, n\}$.

By Corollary \ref{random vertices} we have
\begin{equation}  \label{Rademacher average}
E_r(F) = \sum_{k=1}^{k_1} E_r(F_k)
   \ge c \sum_{k=1}^{k_1} 2^{-k} \sqrt{\log N(k)}
   \ge c k_1
   \ge c_1 \log n.
\end{equation}
To estimate the combinatorial dimension of $F$, fix a $t = 2^{-k}$
with $0 \le k \le k_1$. Then
\[
\vc(F, t)  \le  \vc \Big( \sum_{l=1}^{k+1} F_l, \frac{t}{2} \Big),
\]
because the diameter of $\sum_{l=k+2}^{k_1} F_l$ in the $L_\infty$
norm is at most $t/2$. Obviously, by definition of the
combinatorial dimension, for any finite set $F$ we have $\vc(F,t)
\le \log_2 |F|$, so
\begin{equation} \label{vc example}
\vc(F, t) \le \log_2 \Big| \sum_{l=1}^{k+1} F_l \Big| =
\sum_{l=1}^{k+1} \log_2 |F_l| = \sum_{l=1}^{k+1} 4^l \le C 4^k.
\end{equation}
This shows that $t \sqrt{\vc(F,t)} \le C_1$ for all $t \ge 2^{-k_1}$.
Since $2^{-k_1} \le 2/\sqrt{n}$ and clearly
$t \sqrt{\vc(F,t)} \le t \sqrt{n} \le 2$ for all $t \le 2/\sqrt{n}$,
we conclude that
\begin{equation}                \label{vc small}
\sup_{t > 0} t \sqrt{\vc(F,t)}  \le  C_1.
\end{equation}

Also, since $F$ is uniformly bounded by $1$, we have $\vc(F,t) =
0$ for all $t > 1$. Moreover, for all $f,g \in F$ and all $i \in
\{1,\ldots,n\}$, we have $|f(i)-g(i)| \ge 2^{-k_1}$ whenever $f(i)
\neq g(i)$. Hence, $\vc(F,t) = \vc(F,t_1)$ for all $t \le t_1 =
2^{-k_1}$. Thus,
\[
\int_0^\infty \sqrt{\vc(F,t)} \; dt
\le  C \sum_{k=0}^{k_1} 2^{-k} \sqrt{\vc(F, 2^{-k})}
\le C k_1
\le C \log n.
\]
This, \eqref{Rademacher average} and Theorem \ref{supremum}
imply that the leftmost and the middle quantities in \eqref{eq no sudakov}
are both equivalent to $\log n$ up to an absolute constant factor.
Together with \eqref{vc small}, this completes the proof.
\endproof

\section{Sections of convex bodies}                   \label{s:sections}

First applications of entropy inequalities involving the
combinatorial dimension to geometric functional analysis are due
to M.~Talargand \cite{T 92}. Using his entropy inequality (which
\eqref{eq supremum} strengthens) he proved that for Banach spaces
infratype $p$ implies type $p$ ($1 < p < 2$).
He also proved classical Elton's Theorem with asymptotics
that fell short from optimal, improving
earlier estimates by J.Elton \cite{E} and A.Pajor \cite{Pa};
the optimal asymptotics were found in \cite{MV} using \eqref{Lp}.

Here we will use new covering results to find nice coordinate
sections of a general convex body (for simplicity, we will assume
that the body is symmetric with respect to the origin).
Our main result is related to three classical results in
geometric functional analysis -- Dvoretzky Theorem in the form of
V.~Milman (see \cite{MS} 4.2), Bourgain-Tzafriri's Principle of the
Restricted Invertibility \cite{BT 87} and Elton's Theorem (\cite{E},
see also \cite{Pa}, \cite{T 92}, \cite{MV}).

By $B_p^n$ we denote the unit ball of $l_p^n$, that is the set of
all $x \in \R^n$ such that $\sum_1^n |x(i)|^p \le 1$. Let $K$ be a
convex body symmetric with respet to the origin. Its average size
is measured by $M_K = \int_{S^{n-1}} \|x\|_K \; d\s(x)$,
where $\s$ is the normalized Lebesgue measure on the sphere
$S^{n-1}$ and $\|x\|_K$ denotes the Minkowski functional of $K$
(the seminorm whose unit ball is $K$).

\begin{theorem}[Dvoretzky Theorem, see \cite{MS}]
  Let $K$ be a symmetric convex body in $\R^n$ containing $B_2^n$.
  Then there exists a subspace $E$ in $\R^n$ of dimension
  $k \ge c M_K^2 n$ and such that
  \begin{equation}              \label{eq dvoretzky}
    c(B_2^n \cap E)
    \subseteq M_K (K \cap E)
    \subseteq  C(B_2^n \cap E).
  \end{equation}
  Moreover, a random subspace $E$ taken uniformly in the
  Grassmanian $G_{n,k}$ satisfies \eqref{eq dvoretzky} with
  probability at least $1 - e^{-ck}$.
\end{theorem}

Next theorem, the Principle of the Restricted Invertibility due to
J.~Bourgain and L.~Tzafriri, is the first and probably the most used
result from the extensive paper \cite{BT 87}. By $(e_i)$ we
denote the canonical basis of $\R^n$.

\begin{theorem}[J.~Bourgain and L.~Tzafriri \cite{BT 87}]
  Let $T : l_2^n \to l_2^n$ be a linear operator with $\|T e_i\| \ge 1$
  for all $i$. Then there exists a subset $\sigma$ of
  $\{1, \ldots, n\}$ of size $|\s| \ge cn / \|T\|^2$
  and such that
  $$
  \|Tx\|  \ge  c \|x\| \ \ \
  \text{for all $x \in \R^\s$.}
  $$
\end{theorem}

Denote by $(\e_n)$ Rademacher random variables,
i.e. sequence of independent symmetric $\pm 1$ valued random variables.

\begin{theorem}[J.~Elton \cite{E}]
   Let $x_1, \ldots, x_n$ be vectors in a real Banach space, satisfying
   $$
   \forall i \ \|x_i\| \le 1 \
   \ \ \ \text{and} \ \ \
   \E \Big\| \sum_{i=1}^n \e_i x_i \Big\|  \ge  \d n
   $$
   for some number $\d > 0$.
   Then there exists a subset $\s \subset \{ 1, \ldots, n \}$
   of cardinality $|\s|  \ge  c_1(\d) n$ such that
   $$
   \Big\| \sum_{i \in \s} a_i x_i \Big\|
   \ge   c_2(\d) \sum_{i \in \s} |a_i|
   $$
   for all real numbers $(a_i)$.
\end{theorem}
The best possible asymptotics is known:
$c_1(\d) \asymp \d^2$ and  $c_2(\d) \asymp \d$ \cite{MV}.

Now we state our main result. By $(g_i)$ we denote a sequence of
independent normalized Gaussian random variables.

\begin{theorem}                                           \label{main sec}
   Let $x_1, \ldots, x_n$ be vectors in a real Banach space, satisfying
   \begin{equation}                                  \label{x assumptions}
     \Big\| \sum_{i=1}^n a_i x_i \Big\|
     \le  \sqrt{n} \Big( \sum_{i=1}^n |a_i|^2 \Big)^{1/2}
     \ \ \ \text{and} \ \ \
     \E \Big\| \sum_{i=1}^n g_i x_i \Big\|  \ge  \d n
   \end{equation}
   for all real numbers $(a_i)$ and for some number $\d > 0$.
   Then there exist two numbers $s > 0$ and $c\d \le t \le 1$
   connected by the inequality
   $s t  \ge  c \d / \log^{3/2}(2/\d)$
   and a subset $\s$ of $\{1, \ldots, n\}$ of size $|\s| \ge s^2 n$\
   such that
   \begin{equation}                               \label{on sigma}
     \Big\| \sum_{i \in \s} a_i x_i \Big\|
     \ge   c t \sum_{i \in \s} |a_i|
   \end{equation}
   for all real numbers $(a_i)$.
\end{theorem}

The first assumption in \eqref{x assumptions} is satisfied in
particular if $\|x_i\| \le 1$ $\forall i$.
Also, since $t \le 1$, we always have
$s \ge c \d / \log^{3/2}(2/\d)$.
This instantly recovers Elton's Theorem.

Next, Theorem \ref{main sec} essentially extends the
Bourgan-Tzafriri's principle of restricted invertibility
to operators $T : l_2^n \to X$ acting into arbitrary Banach
space $X$. The average size of $T$ is measured by
its $\ell$-norm defined as $\ell(T)^2 = \E \|Tg\|^2$,
where $g = (g_1, \ldots, g_n)$. If $X$ is a Hilbert
space, then $\ell(T)$ equals the Hilbert-Schmidt norm of $T$.

\begin{corollary}[General Principle of the Restricted Invertibility]
 \label{gen PRI}
  Let $T : l_2^n \to X$ be a linear operator with $\ell(T) \ge \sqrt{n}$,
  where $X$ is a Banach space. Let $\a = c \log^{-3/2} (2\|T\|)$.
  Then there exists a subset $\sigma$ of $\{1, \ldots, n\}$
  of size $|\s| \ge c \a^2 n / \|T\|^2$ and such that
  $$
  \|Tx\|  \ge  \a |\s|^{-1/2} \|x\|_{l_1^\s} \ \ \
  \text{for all $x \in \R^\s$.}
  $$
\end{corollary}

If $X$ is a Hilbert space, the condition $\ell(T) \ge \sqrt{n}$ is
satisfied, for example, if $\|x_i\| \ge 1$ for all $i$. In that
case $ |\s|^{-1/2} \|x\|_{l_1^\s}$ in the conclusion can be
improved to $|x\| = \|x\|_{l_2^n}$ via the Grothendieck
factorization (we will do this below). This recovers
Bourgain-Tzarfiri's Theorem up to the logarithmic factor $\a$.

\qquad

\noindent {\bf Proof of Corollary \ref{gen PRI}. } We apply
Theorem \ref{main sec} to the vectors $x_i =
\frac{\sqrt{n}}{\|T\|} \, T e_i$, $i = 1, \ldots, n$. Then for any
$a_1, \ldots a_n$
$$
\left \| \sum_{i=1}^n a_i x_i \right \| =  \frac{ \sqrt{n}}{\|T\|}
\cdot \left \|T( \sum_{i=1}^n a_i e_i) \right \| \le  \sqrt{n}
\cdot \|(a_1, \ldots a_n)\|_{\ell_2^n}.
$$
Since by Kahane's inequality $\ell(T) \le C \, \E \|Tg\|$, the
second assumption in \eqref{x assumptions} holds with $\d =
c/\|T\|$. Hence there exist numbers $c/\|T\| \le t \le 1$ and $s$
satisfying
$$
s t  \ge  c \d / \log^{3/2}(2/\d)
$$
and a subset $\sigma$ of $\{1, \ldots, n\}$ of size $|\s| \ge s^2
n$ so that we have (multiplying both sides by $|\s|^{-1/2}$)
$$
\frac{C\sqrt{n/|\s|}}{\|T\|}
  \Big\| \sum_{i \in \s} a_i Te_i \Big\|
\ge   c t |\s|^{-1/2} \sum_{i \in \s} |a_i|
\ \ \ \text{for all real numbers $(a_i)$.}
$$
Since $t \le 1$, we have
$s \ge c \d / \log^{3/2}(2/\d)  \ge  \a / \|T\|$ and
consequently $|\s| \ge \a^2 n / \|T\|^2$ as required. As
$|\s| \ge s^2 n$, we have $s \le \sqrt{|\s|/n}$,
hence
$$
\frac{\sqrt{n/|\s|}}{\|T\|}
\le  \frac{1}{\|T\| s}
=    \d / s
\le  c^{-1} t \log^{3/2}(2/\d)
= t / \a.
$$
Hence
$$
\Big\| \sum_{i \in \s} a_i Te_i \Big\|
\ge  \a |\s|^{-1/2} \sum_{i \in \s} |a_i|
\ \ \ \text{for all real numbers $(a_i)$}
$$
as required.
\endproof

To get the actual invertibility of $T : l_2^n \to X$
one can use the Grothendieck factorization.
Remarkably, this step works not only for $X$ being a Hilbert space
but for a much larger class of spaces, namely for those of type
$2$.

\begin{definition}
  A Banach space $X$ has {\em type $2$}
  if there exists a constant $M$ such that the inequality
  $$
  \E \Big\| \sum \e_i x_i \Big\|
  \le M \Big( \sum \|x_i\|^2 \Big)^{1/2}
  $$
  holds for all finite sequences of vectors $(x_i)$ in $X$.
  The minimal possible constant $M$ is called the {\em type $2$ constant of $X$}
  and is denoted by $T_2(X)$.
\end{definition}

An important example of spaces that have type $2$ are
all $L_p$-spaces ($2 \le p < \infty$) and their subspaces.

\begin{lemma}[Grothendieck Factorization, see e.g. \cite{LT} 15.4]
   Let $S : E \to \R^m$ be a linear mapping, where $E$ is a Banach
   space
   of type $2$. Then there exists a subset
   $\eta$ of $\{1, \ldots, n\}$
   of size $|\eta| \ge m/2$ and such that
  $$
  \|P_\eta S\|_{E \to l_2^m}
  \le C \, T_2(X) \, m^{-1/2} \|S\|_{E \to l_1^m}
  $$
  where $P_\eta$ is the coordinate projection in $\R^m$ onto $\R^\eta$.
\end{lemma}

Applying this lemma to the inverse of $T$ on its range, we obtain

\begin{corollary}[Restricted Invertibility under type $2$]
  Let $T : l_2^n \to X$ be a linear operator with $\ell(T) \ge \sqrt{n}$,
  where $X$ is a Banach space of type $2$.
  Let $\a = c \log^{-3/2} (2\|T\|)$.
  Then there exists a subset $\sigma$ of $\{1, \ldots, n\}$
  of size $|\s| \ge \a^2 n / \|T\|^2$ and such that
  $$
  \|Tx\|  \ge  \a \, T_2(X)^{-1} \|x\| \ \ \
  \text{for all $x \in \R^\s$.}
  $$
\end{corollary}
For $X = l_2^n$ this recovers Bourgain-Tzafriri Theorem up to the
logarithmic factor $\a$.

\qquad

\proof By Corollary \ref{gen PRI}, the operator $T$ is
invertible on the subspace $E = T(\R^\s)$ of $X$, and its inverse
$S = T^{-1} : E \to \R^\s$ has norm $\|S\|_{E \to l_1^\s}  \le
\a^{-1} |\s|^{1/2}$. By the Grothendieck factorization, we find
a subset $\eta \subset \s$ of size $|\eta| \ge \frac{1}{2} |\s|$
and such that
$$
\|P_\eta S\|_{E \to l_2^\eta}  \le  C \a^{-1} T_2(X).
$$
This means that $\|Tx\|_X  \ge  c \a \, T_2(X)^{-1} \|x\|_{l_2^n}$
for all $x \in \eta$.
\endproof

\qquad

Finally, the general Principle of the Restricted Invertibility
rewritten in geomteric terms gives a result related to
Dvoretzky Theorem.

\begin{corollary}                                  \label{coord Dvoretzky}
  Let $K$ be a symmetric convex body in $\R^n$ containing $B_2^n$.
  Let $M = M_K \log^{-3/2} (2/M_K)$.
  Then there exists a subset $\s$ of $\{1, \ldots, n\}$
  of size $|\s|  \ge  c M^2 n$ and such that
  \begin{equation}                                  \label{section in B1}
    c M  (K \cap \R^\s)  \subset  \sqrt{|\s|} B_1^\s.
  \end{equation}
\end{corollary}
Here $B_1^\s$ denotes the unit ball of $l_1^\s$.

\qquad

\proof We apply the general Principle of Restricted Invertibility
in the space $X = (\R^n, \|\cdot\|_K)$. We have $\ell(\id : l_2^n
\to X) = \sqrt{n} M_K$ (see \cite{TJ} (12.7)) and $\|\id : l_2^n
\to X\| \le 1$ because $K$ contains $B_2^n$. Hence for the
operator $T = c (M_K)^{-1}\id : l_2^n \to X$ we have $\ell(T) \ge
\sqrt{n}$ and $\|T\| \le C/M_K$. The application of Corollary
\ref{gen PRI} completes the proof.
\endproof

\qquad

A link to Dvoretzky Theorem is provided by  a result of Kashin
\cite{K 77} (see also \cite{S}) that the cross-polytope $\sqrt{k}
B_1^k$ has a Euclidean section of proportional dimension.
Precisely, there exists a subspace $E$ in $\R^k$ of dimension at
least $k/2$ and such that
$$
(B_2^k \cap E) \subseteq  (\sqrt{k} B_1^k \cap E) \subseteq  C
(B_2^k \cap E).
$$
Actually, a random subspace $E$ taken uniformly in the Grassmanian
satisfies this with probability at least $1 - e^{-ck}$.

Taking such  random section of both sides of \eqref{section in B1}
we get $M(K \cap E) \subseteq  C (B_2^n \cap E)$, which recovers
the second inclusion in Dvoretzky Theorem up to a logarithmic
factor. The novelty of \eqref{section in B1} is that the section
is coordinate. This might be important for future applications.

\qquad

\remark
Corollary \ref{coord Dvoretzky} may fail for any set of size $|\s|
\asymp  M^2 n$, even though it must hold for some larger set.
Indeed, for $K = a \sqrt{n} B_1^n$ with some large parameter $a$
we have $M \sim a^{-1} \log^{-3/2} a$. Any set $\s$ for which
Corollary \ref{coord Dvoretzky} holds satisfies $\sqrt{|\s|}
B_1^\s \supseteq  c a^{-1} (\log^{-3/2} a) (a \sqrt{n} B_1^n \cap
\R^\s) = c (\log^{-3/2} a) \sqrt{n} B_1^\s$, so $|\s| \ge
(\log^{-3/2} a) n$. This is much larger than $M^2 n  \asymp  a^{-2}
(\log^{-3} a) n$. In particular, Corollary \ref{coord Dvoretzky}
fails for any set of size $|\s| \sim  M^2 n$.

\qquad

\noindent {\bf Proof of Theorem \ref{main sec}. } By a slight
perturbation we may assume that the vectors $x_i$ are linearly
independent, and by applying appropriate linear transformation we may
further assume that $X = (\R^n, \|\cdot\|_K)$ where $K$ is a
symmetric convex body in $\R^n$ and that $x_i = e_i$, the
canonical vector basis in $\R^n$. We then rewrite the assumptions
as $\frac{1}{\sqrt{n}} B_2^n  \subseteq K$, $\E \|g\|_K \ge \d n$.
Then for the polar body
$A = K^\circ = \{x \in \R^n : \< x,y \> \le 1 \ \forall y \in K \}$
we have
\begin{equation}                             \label{on A}
  A \subseteq \sqrt{n} B_2^n, \ \ \
  E := \E \sup_{x \in A} \sum_{i=1}^n g_i x(i)  \ge  \d n.
\end{equation}

Although Theorem \ref{supremum} can be used to estimate $E$,
we will need to have some control on the upper limit
in the integral \eqref{eq supremum}. This can be done as follows.
Noting that $E(A) = n^{-1/2} E$,
we bound the expectation in \eqref{on A} by Dudley's entropy inequality
\eqref{dudley}:
\begin{equation}                \label{dudley worked}
E  \le  C \sqrt{n} \int_{cE/n}^1 \sqrt{\log N(A, tD)} \; dt
\end{equation}
where $D = \Ball(L_2^n) = \sqrt{n} B_2^n$.
The upper limit in the integral is $1$ because $A \subseteq D$, so the
integrand vanishes for $t > 1$.
By Theorem \ref{by ellipsoids},
$$
N(A, tD) = N(t^{-1} A, D) \le  \Sigma(C t^{-1} A)^2.
$$
Since $C t^{-1} A \subseteq C t^{-1} \cdot \Ball(L_2^n)$,
Lemma \ref{A in L1} gives
$$
\Sigma(C t^{-1} A)  \le  \Big( \frac{Cn}{t \, v(t)} \Big)^{C v(t)},
$$
where $v(t) = \vc(C t^{-1} A)$. Hence
$$
\log N(A, tD) \le  C v(t) \log \Big( \frac{C n}{t \, v(t)} \Big).
$$
Using this in Dudley's entropy inequality \eqref{dudley worked},
we obtain
$$
E  \le  C \sqrt{n} \int_{cE/n}^1
        \sqrt{ v(t) \log \Big( \frac{Cn}{t \, v(t)} \Big) }.
        \; dt
$$
Let $s(t)^2 = v(t)/n$. Since $s(t) \le 1$ and $E \ge \d n$, we
have
$$
c \d  \le  \int_{c \d}^1
  s(t) \sqrt{ \log \Big (\frac{1}{t \, s(t)} \Big) } \; dt.
$$
Comparing the integrand to that of
$$
\log(1 / c \d) = \int_{c \d}^1  \frac{1}{t} \; dt
$$
we conclude that there exists a number $c \d \le t \le 1$ such
that
$$
s(t) \sqrt{ \log \Big (\frac{1}{t \, s(t)} \Big) } \ge  \frac{c
\d}{t \log(1 / c \d)}.
$$
Multiplying both sides by $t$, we obtain
$$
t \, s(t) \ge  \frac{c\d}{\log(1/c\d)}
     \Big/ \sqrt{ \log \Big( \frac{\log(1/c\d)}{c\d} \Big) }
\ge  \frac{c \d}{\log^{3/2}(2/\d)}.
$$

It remains to interpret $v(t)$.
By the symmetry of $A$, $v(t)$ is the maximal rank of a coordinate
projection $P$ in $\R^n$ such that
$P(C t^{-1} A)  \supseteq  P( \frac{1}{2} B_\infty^n)$.
Let $\R^\s$ be the range of $P$; then $|\s| = v(t) = s(t)^2 n$.
By duality, the inclusion above is equivalent to
$C^{-1} t K \cap \R^\s  \subseteq  2 B_1^n$.
Equivalently, $\|x\|_K  \ge  C^{-1} t \|x\|_{l_1^n}$
for all $x \in E$. This is precisely the conclusion \eqref{on sigma}.
The proof is complete.
\endproof

\remark Although the first assumption in \eqref{x assumptions} is
rather nonrestrictive, it can further be weakened. Tracing where it
was used in the proof (in Lemma \ref{A in L1}) we see that only
"average" volumetric properties of $K$ matter. We leave details to
the interested reader.

\end{document}